\documentclass{amsart}
\usepackage{amsmath,amssymb,amsxtra,amsthm,amscd}
\usepackage{lmodern}
\usepackage[T1]{fontenc}
\usepackage[all,cmtip]{xy}
\usepackage{rotating}
\usepackage{microtype}
\usepackage[linecolor=yellow,color=yellow,bordercolor=yellow]{todonotes}
\newif\ifShowLabels
\ShowLabelstrue
\newcommand{\TeXref}[1]{
\marginpar{\scriptsize \texttt{#1}}}

\DeclareMathOperator{\B}{\mathbf{B}}

         \newcommand{\BGamma}{\B_{\Gamma}}

         \newcommand{\shlf}{\sideset{^{s}}{^{\textit{lf}}}\h}

\DeclareMathOperator{\EGamma}{\boldsymbol{E}\boldsymbol{\Gamma}}

\DeclareMathOperator{\fil}{fil}

\DeclareMathOperator{\Free}{\mathbf{Free}}

\DeclareMathOperator{\Fun}{Fun}

\DeclareMathOperator{\h}{\mathit{h}}
\DeclareMathOperator{\hlf}{\mathit{h}^{\textit{lf}}}

\DeclareMathOperator{\holimplain}{holim}
\DeclareMathOperator{\Hom}{Hom}

\DeclareMathOperator{\id}{id}

\DeclareMathOperator{\K}{\mathit{K}}
         \newcommand{\Knc}{\K^{-\infty}}

\DeclareMathOperator{\LI}{\mathbf{L}}
\DeclareMathOperator{\SI}{\mathbf{S}}

\DeclareMathOperator{\BLI}{\mathbf{BL}}
\DeclareMathOperator{\BSI}{\mathbf{BS}}

\DeclareMathOperator{\Mod}{\mathbf{Mod}}

\DeclareMathOperator{\Or}{Or}

\DeclareMathOperator{\point}{point}

\DeclareMathOperator{\U}{\mathbf{U}}

\DeclareMathOperator{\W}{\mathbf{W}}

\DeclareMathOperator*{\one}{1}
\newcommand{\onehatplace}[1]%           % "frown" over 1
{ \one^{\substack{#1 \\ \frown}} }

\DeclareMathOperator*{\bones}{\times}
\newcommand{\undertimes}[1]%            % X with stuff under it
{ \bones_{#1} }

\DeclareMathOperator*{\bowl}{\cup}
\newcommand{\undercup}[1]%              % small U with stuff under it
{ \bowl_{#1} }

\DeclareMathOperator*{\arch}{\cap}
\newcommand{\undercap}[1]%              % small "cap" with stuff under it
{ \arch_{#1} }

\newcommand{\pull}%                              % extension for arrows
{\!\!\! -\!\!\! -\!\!\! -\!\!\!}

\DeclareMathOperator*{\holimprep}{holim}                       % holim
\newcommand{\holim}[1]%
{\displaystyle\holimprep_{\substack{\leftarrow \pull - \\ #1}} \, }

\DeclareMathOperator*{\hocolimprep}{hocolim}                   % hocolim
\newcommand{\hocolim}[1]%
{\displaystyle\hocolimprep_{\substack{- \pull \rightarrow \\ #1}} \, }

\DeclareMathOperator*{\plainlim}{lim}                           % lim
\newcommand{\contralim}[1]%
{\displaystyle\plainlim_{\substack{\leftarrow \pull - \\ #1}} \, }

\DeclareMathOperator*{\plaincolim}{colim}                       % colim
\newcommand{\colim}[1]%
{\displaystyle\plaincolim_{\substack{- \pull \rightarrow \\ #1}} \, }

\DeclareMathOperator*{\laxlimplain}{laxlim}                     % laxlim
\newcommand{\laxlim}[1]%
{\displaystyle\laxlimplain_{\substack{\leftarrow \pull - \\ #1}} \, }

\providecommand{\bysame}{\makebox[3em]{\hrulefill}\thinspace}

\providecommand{\cal}{\it}     % substitute for \cal

\theoremstyle{plain}
\newtheorem{Thm}{Theorem}[section]

\newtheorem{Cor}[Thm]{Corollary}

\newtheorem{Lem}[Thm]{Lemma}
\newtheorem{Prop}[Thm]{Proposition}

\theoremstyle{definition}
\newtheorem{Def}[Thm]{Definition}

\newtheorem{Ex}[Thm]{Example}

\newtheorem{Rem}[Thm]{Remark}

\newtheorem*{Sum}{Summary}
\newtheorem*{Cav}{Caveat}

\theoremstyle{remark}
\newtheorem{Not}[Thm]{Notation}

\newtheoremstyle{freestylethm}{6pt}{6pt}{\itshape}{}%
                {\bfseries}{}{.5em}{\thmnote{#3}}
\theoremstyle{freestylethm}

\newcommand{\SecRef}[2]{\section{#1}\label{S:#2}%
\ifShowLabels \TeXref{{S:#2}} \fi}
\newcommand{\SSecRef}[2]{\subsection{#1}\label{SS:#2}%
\ifShowLabels \TeXref{{SS:#2}} \fi}

\newcommand{\refS}[1]{\textup{\ref{S:#1}}}
\newcommand{\refSS}[1]{\textup{\ref{SS:#1}}}

\newcommand{\refT}[1]{\textup{\ref{T:#1}}}

\newcommand{\refD}[1]{\textup{\ref{D:#1}}}
\newcommand{\refC}[1]{\textup{\ref{C:#1}}}

\newenvironment{ThmRef}[1]%
{ \begin{Thm} \label{T:#1}
\ifShowLabels \TeXref{T:#1} \fi }%
{ \end{Thm} }
\newenvironment{DefRef}[1]%
{ \begin{Def} \label{D:#1}
\ifShowLabels \TeXref{D:#1} \fi }%
{ \end{Def} }
{ \begin{Lem} \label{L:#1}
\ifShowLabels \TeXref{L:#1} \fi }%
{ \end{Lem} }
\newenvironment{CorRef}[1]%
{ \begin{Cor} \label{C:#1}
\ifShowLabels \TeXref{C:#1} \fi }%
{ \end{Cor} }
\newenvironment{RemRef}[1]%
{ \begin{Rem} \label{R:#1}
\ifShowLabels \TeXref{R:#1} \fi }%
{ \end{Rem} }
\newenvironment{PropRef}[1]%
{ \begin{Prop} \label{P:#1}
\ifShowLabels \TeXref{P:#1} \fi }%
{ \end{Prop} }
{ \begin{Ex} \label{E:#1}
\ifShowLabels \TeXref{E:#1} \fi  }%
{ \end{Ex} }
\newenvironment{NotRef}[1]%
{ \begin{Not} \label{N:#1}
\ifShowLabels \TeXref{N:#1} \fi }%
{ \end{Not} }

\newenvironment{ThmRefName}[2]%
{ \begin{Thm} [#2]\label{T:#1}
\ifShowLabels \TeXref{T:#1} \fi }%
{ \end{Thm} }
\newenvironment{DefRefName}[2]%
{ \begin{Def} [#2]\label{D:#1}
\ifShowLabels \TeXref{D:#1} \fi }%
{ \end{Def} }
{ \begin{Lem} [#2]\label{L:#1}
\ifShowLabels \TeXref{L:#1} \fi }%
{ \end{Lem} }
{ \begin{Cor} [#2]\label{C:#1}
\ifShowLabels \TeXref{C:#1} \fi }%
{ \end{Cor} }
{ \begin{Rem} [#2]\label{R:#1}
\ifShowLabels \TeXref{R:#1} \fi }%
{ \end{Rem} }
{ \begin{Prop} [#2]\label{P:#1}
\ifShowLabels \TeXref{P:#1} \fi }%
{ \end{Prop} }
\newenvironment{ExRefName}[2]%
{ \begin{Ex} [#2]\label{E:#1}
\ifShowLabels \TeXref{E:#1} \fi }%
{ \end{Ex} }

 \ShowLabelsfalse

\begin{document}

\title[Equivariant stable homotopy methods]{Equivariant stable homotopy methods in the\\ algebraic K-theory of infinite groups}
\author{Gunnar Carlsson}
\address{Department of Mathematics\\ Stanford University\\ Stanford\\ CA 94305}
\email{gunnar@math.stanford.edu}
\author[Boris Goldfarb]{Boris Goldfarb}
\address{Department of Mathematics and Statistics\\ SUNY\\ Albany\\ NY 12222}
\email{goldfarb@math.albany.edu}
\date{January 4, 2019}

\begin{abstract}
Equivariant homotopy methods developed over the last 20 years lead to recent breakthroughs in the Borel isomorphism conjectures for Loday assembly maps in K- and L-theories. 
An important consequence of these algebraic conjectures is the topological rigidity of compact aspherical manifolds.
Our goal is to strip the basic idea to the core and follow the evolution over time in order to explain the advantages of the flexible state that exists today.
We end with an outline of the proof of the Borel conjecture in algebraic K-theory for groups of finite asymptotic dimension. 
\end{abstract}

\maketitle

\tableofcontents

\SecRef{Introduction}{intro}

We start with a very general idea in algebraic topology.
Suppose we are given a continuous map of spaces
$f \colon \underline{S} \to \underline{T}$.
One is generally interested in the homomorphism that is induced on the homotopy groups 
$f_n \colon \pi_n (\underline{S}) \to \pi_n (\underline{T})$.
The basic question is whether $f_n$ is an isomorphism for a given integer $n$, or an injection, or a surjection.
If $f$ induces isomorphisms $f_n$ for all $n$,
we refer to $f$ as a weak homotopy equivalence or simply an equivalence.

Since we will be working on the level of spaces, in order to translate such properties 
via the functor $\pi_n$ we will want to establish split injectivity and split surjectivity. 
This means, in the instance of split injectivity, that there should be another map $g \colon \underline{T} \to \underline{S}$ called a \textit{splitting} such that $g \circ f \simeq \id$ on $\underline{S}$.
In this case, for any homotopy invariant functor $\Phi$ on the category of spaces, $\Phi (g) \circ \Phi (f) = \id$ on $\Phi (S)$, and so $\Phi (f)$ is also a split injection.
This can be applied to the functors $\Phi = \pi_n$.
Also, attempting to do this on the level of spaces means we are shooting for split injectivity for all integers $n$.

The basic way to establish that $f$ is a split injection is thus to \textit{engineer} another spectrum $\underline{M}$ which receives two maps $s \colon \underline{S} \to \underline{M}$ and $t \colon \underline{T} \to \underline{M}$ which form a commutative traingle
\[
\xymatrix{
 \underline{S}   \ar[ddr]^-{s} \ar[rr]^-{f}
&&\underline{T}  \ar[ddl]_-{t} \\
\\
& \underline{M}
}
\]
while $s$ is also an equivalence.

In order to address split surjectivity of $f$, one could simply dualize this template.  
Another method that we will eventually use in this paper can be applied to a map that has just been shown to be a split injection in the triangle above.
One may attempt to use the same idea on the splitting map $t$:
if $t$ is split injective then it is an equivalence, and that makes $f$ an equivalence.

In many areas of mathematics, the crucial invariants are related to stable homotopy groups.  These groups have been organized into
so-called generalized homology and cohomology theories represented as additive theories by objects called connective spectra. For a good background in the theory of spectra the reader is advised to consult any of the books 
\cite{fA:78,fA:95,gCjM:95,jpM:77}, or \cite[Chapter II]{yR:98}.
The 0-level of a spectrum is a topological space with an infinite loop space structure. There is in fact an equivalence between spectra and the infinite loop spaces. Everything said so far can be stated in the context of spectra and so be applied to stable homotopy group functors on the category of spectra.  Most of the paper can be understood if the term \textit{spectrum} and maps between spectra are replaced with \textit{space} or, specifically, infinite loop space and maps between them.

The equivariant homotopy theoretic angle arises from the fact that very often $f$ can be modeled, up to homotopy, as a fixed point map.
So suppose there are left actions of a group $\Gamma$ on two spectra $S$ and $T$ and the map $\phi \colon S \to T$ which is $\Gamma$-equivariant.
Then one has the induced map on the fixed point spectra 
$\phi^\Gamma \colon S^\Gamma \to T^\Gamma$.
We are assuming that the choices are such that there are weak equivalences $\underline{S} = S^\Gamma$, $\underline{T} = T^\Gamma$, and $f = \phi^\Gamma$ up to homotopy.

It is important to point out that properties of $\phi$ such as being an equivalence or split injection do not in principle get inherited by $f$.
What may be initially surprising is that selecting a good target $\underline{M}$ gets easier in the equivariant setting.
Here one would naturally look for an equivariant spectrum $M$ and for equivariant maps $\sigma \colon S \to M$ and $\tau \colon T \to M$ which fit a commutative triangle
\[
\xymatrix{
 {S}   \ar[ddr]^-{\sigma} \ar[rr]^-{\phi}
&&{T}  \ar[ddl]_-{\tau} \\
\\
& {M}
}
\]
and such that the fixed point map $\sigma^\Gamma$ is an equivalence.
\textit{It is a phenomenon in equivariant topology that such choices are easier to make than guessing 
$\underline{M}$ directly.}

\medskip 

This paper will follow the historical path from the basic idea of using homotopy fixed points
to prove injectivity of the assembly (the Novikov Conjecture) through further refinements to the recent proof of surjectivity (the Borel Conjecture) in all of the cases where the Novikov Conjecture is currently known.
One word of warning: we concentrate on the development of technology rather than a comprehensive overview of the accomplishments in the subject of assembly maps.  
However, in the last section, we will explore connections to the recent work on the Farrell-Jones conjecture.

\SecRef{Bounded K-theory and the Loday assembly map}{bddkth}

\SSecRef{Homotopy fixed point method}{hfp}

The application of homotopy fixed points to splitting the assembly map was first realized in \cite{gC:95} and was used extensively since then by many authors.

Since we use the so-called naive equivariant spectra, the reader will not need a background in modern equivariant stable homotopy theory.  We advise to consult the following references \cite{gC:93,gC:05} for some explicit elementary details as needed by the interested reader.

For a spectrum $S$ with a left $\Gamma$-action, the function spectrum $\Fun (E\Gamma, S)$ has the induced left $\Gamma$-action.  Here $E\Gamma$ stands for the universal contractible free $\Gamma$-space.  The left action is given by $g \cdot \rho = g \circ \rho \circ g^{-1}$. 
The \textit{homotopy fixed points} of $S$, denoted as $S^{h\Gamma}$, are defined as the fixed points $\Fun (E\Gamma, S)^\Gamma$.
They are simply the equivariant maps from $E\Gamma$ to $S$, the maps that respect the actions in the domain and the codomain in the evident way.
It is a general fact that the induced map $\phi^{h\Gamma} = \Fun (E\Gamma, \phi)^\Gamma \colon S^{h\Gamma} \to T^{h\Gamma}$ is an equivalence whenever 
$\phi$ is a (nonequivariant) equivalence. 
This is in contrast to the fact that $\phi$ being an equivalence has no bearing on the properties of the fixed point map $\phi^{\Gamma}$.

There is a canonical map that relates the fixed points and the homotopy fixed points.
Notice that the fixed point spectrum of $S$ can be viewed along the same lines as the equivariant maps $\Fun (\point, S)^\Gamma$.
The map $r_S \colon S^\Gamma \to S^{h\Gamma}$ is induced by the (equivariant) collapse $E\Gamma \to \point$ giving $\Fun (\point, S)^\Gamma \to \Fun (E\Gamma, S)^\Gamma$.
The problem of analyzing $r_S$ is generally a fundamental and difficult problem.
The advantages are well-known, see for example the discussion around Proposition 1.8 in \cite{gC:92}.
They were generally promoted by Thomason \cite{rT:83}.

The following is a template for using the homotopy fixed points in order to split the fixed point map.
For an equivariant $\phi$ there is always a commutative square
\[
\xymatrix{
 S^\Gamma   \ar[dd]^-{r_S} \ar[rr]^-{f = \phi^{\Gamma}}
&&T^\Gamma  \ar[dd]^-{r_T} \\
\\
 S^{h\Gamma} \ar[rr]^-{\phi^{h\Gamma}}
&&T^{h\Gamma}
}
\]
Let us assume that (1) $r_S$ happens to be an equivalence.  This, combined with the fact that (2) the equivariant map $\phi$ is nonequivariantly an equivalence, gives that the composition $\phi^{h\Gamma} \circ r_S$ is a composition of two equivalences and is, therefore, an equivalence.
By general nonsense, the other composition $r_T \circ f$ should have the first map a split injection and the second map a split surjection.

\SSecRef{The Loday assembly map}{loday}

It is time to look at Loday's assembly map \cite{jlL:76} in algebraic K-theory and its importance in geometric topology.
We will then describe in section \refSS{BBB} how to model the assembly as a fixed point map.
 
We would like the material in this paper to be accessible to someone with a cursory background in K-theory and willing to accept some high-level statements on faith.  It should be enough to believe that K-theory is a functor that produces a spectrum for a ring, or more generally an exact category, and to have seen very classical treatments of the low-dimensional groups $K_0$ and $K_1$.  For example, several chapters from the textbook \cite{jR:94} should suffice, and the book \cite{cW:13} is certainly sufficient.

One can view algebraic K-theory as a functor that produces a spectrum $K(\mathcal{A})$ for a small additive category $\mathcal{A}$. 
For important choices of $\mathcal{A}$ there are meaningful K-groups of such $\mathcal{A}$ which are indexed with negative integers.  Luckily, it is also technically easier to work with the Gersten-Wagoner nonconnective spectrum $K(\mathcal{A})$ whose stable homotopy groups are the negative and positive K-groups of $\mathcal{A}$.
We assume that $K(\mathcal{A})$ is that nonconnective spectrum in the rest of the paper.
For any ring $A$, the K-theory of $A$ is the nonconnective spectrum $K(\mathcal{A})$ where $\mathcal{A}$ is the category of all finitely generated free modules over $A$ with the $\oplus$ operation given by the direct sum.  The most interesting case of K-theory for topologists is this functor applied to the a group ring $R\Gamma$ where $R$ is a ring and $\Gamma$ is a group.  In topological applications, $R$ is usually the integers $\mathbb{Z}$ and the group $\Gamma$ is the fundamental group of a CW-complex.

The assembly map is designed to probe the K-theory spectrum $K(R\Gamma)$ by comparing it to a computationally reasonable group homology spectrum $H(\Gamma, KR)$.
Given a model for the classifying space $B\Gamma$, the latter spectrum is simply the product $B\Gamma_{+} \wedge KR$, where the substript $+$ indicates the base point disjoint from $B\Gamma$.  

Let $\imath \colon  \Gamma \rightarrow GL_1 (R\Gamma)$
be the inclusion of $\Gamma$ in
$(R\Gamma)^{\times} = GL_1(R\Gamma)$.
Then there is a map
\[
\Gamma \times GL_n (R)
\xrightarrow{\imath \times id}
GL_1(R\Gamma) \times GL_n (R)
\xrightarrow{\ \otimes \ }
GL_n(R\Gamma)
\]
defined by
\[
g,(a_{ij}) \longmapsto (g \cdot a_{ij}).
\]
One can apply the classifying space functor $B$,
pass to the limit as $n\rightarrow\infty$,
and apply Quillen's plus construction to induce the map
\[
B\Gamma_+ \wedge BGL (R)^+
\xrightarrow{B\imath^+ \wedge id}
BGL(R\Gamma)^+ \wedge BGL (R)^+
\xrightarrow{\ \gamma \ }
BGL(R\Gamma)^+.
\]
This product is compatible with the infinite loop space
structure of $BGL(\underline {\phantom {S}})^+$.
Delooping of this map results
in the assembly map of spectra
\[
\alpha_K \colon  B\Gamma_+ \wedge K(R) \longrightarrow K(R\Gamma),
\]
where $K(R)$ is the Gersten-Wagoner
non-connective K-theory spectrum of $R$.
This is the {\it Loday assembly map in
algebraic K-theory\/}.

\SSecRef{Geometric applications of the assembly map}{APPL}

There is a good number of reasons
why the study of this map is of importance in geometric topology.
Some of them depend on importance of K-theory computations and use the fact that homology is a computable topological invariant, and the K-theory in contrast is difficult to compute.  So properties of the assembly map, especially in cases when it is an equivalence, allow to use the domain which is homology to derive some information about the target, the K-theory of the group ring.  Other applications depend on the genuine fact that it is the assembly map that is an equivalence.

The first reason we mention is as a way to compute the classical low-dimensional K-groups of ${\Bbb Z} \Gamma$
which contain obstructions to geometric constructions in $M$.
 
Probably the most well-known of these obstructions is the Wall finiteness obstruction.  It is an element in the reduced $K_0$ of ${\Bbb Z} \Gamma$, and if $\alpha_{K,0}$ is an isomorphism then it lifts to the homology, and it can be seen there that the obstruction vanishes.  This happens if $\Gamma$ is the fundamental group of a compact ANR as shown by West \cite{jW:77}.  For the background, general framework, and these kind of computations we refer the interested reader to \cite{kV:89,sW:94,sFaR:01,iH:02}. 

Of course, the Whitehead group, which is the cokernel of $\alpha_{K,1}$ is a well-known home for obstructions such as the Whitehead torsion.  The vanishing of this specific element associated to a homotopy equivalence between compact CW-complexes is a consequence of surjectivity of $\alpha_{K,1}$.  More refined results, such as the fact that the torsion of a homeomorphism has  to vanish, use the assembly even without vanishing of the whole Whitehead group, see \cite{tC:74,aRmY:95,mW:02}.  The paper of Ferry \cite{sF:81} and the survey by Ferry and Ranicki \cite{sFaR:01} give relations to the last paragraph.  

In the category of topological manifolds $M$, the Whitehead group for $\Gamma = \pi_1  M$ contains a number of famous obstructions in addition to the Whitehead torsion. One example is the Farrell-Siebenmann obstruction to fibering a high-dimensional ($\ge 6$) manifold over the circle $S^1$.
There is also the essential involvement of $K({\Bbb Z} \Gamma)$, and the assembly in particular,
in the description of the space of automorphisms
of the manifold. 
Details of this application can be found in \cite{fF:02,fFlJ:91,mWbW:01}.  
 
There is an observation of Waldhausen in his seminal paper \cite[page 139]{fW:78}. He points out that when the assembly map $\alpha_K$ is an equivalence for both of two groups $\Gamma_1$ and $\Gamma_2$, if a homomorphism of groups $\Gamma_1 \to \Gamma_2$ induces an isomorphism on integral homology then their K-groups are canonically isomorphic.  That's what happens when $\Gamma_1$ is a classical knot group, $\Gamma_2$ is the free cyclic group, and the homomorphism is the abelianization homomorphism.

Yet another reason one wants to compute the assembly is a connection with the Novikov and Borel conjectures.
There is a companion construction of the assembly map in algebraic L-theory
\[
\alpha_L \colon  B\Gamma_+ \wedge L(R) \longrightarrow L(R\Gamma),
\]
It turns out that even injectivity of this map has many exciting consequences.
It is known that the classical version of
the Novikov conjecture, the conjecture about higher signatures that
\begin{equation}
f_{\ast} \bigl( L(M)\cap\lbrack M \rbrack \bigr) =
f_{\ast}g_{\ast} \bigl( L(M')\cap\lbrack M' \rbrack \bigr) \in
H_{\ast}(B\Gamma;\Bbb Q)  \tag{$\ddag$}
\end{equation}
whenever $g\colon M'\rightarrow M$ is a homotopy equivalence,
follows from the splitting of the rational assembly map
$\alpha$ in L-theory. 
The assembly naturally maps the rational group homology containing the
signature to the surgery L-group where the image is a priori
homotopy invariant.

The geometric L-theory assembly map was defined by Frank Quinn \cite{fQ:70}.  The algebraic L-theory assembly is due to Andrew Ranicki \cite{aR:80}.
Omitted details can be found in the original papers, or in Ranicki's surveys \cite{aR:95,aR:01}, or the comprehensive monographs \cite{aR:92,sW:94}. 
These references are also a good source for the origin of and motivation for formulating the Novikov conjecture in terms of the assembly in algebraic L-theory.  

The Novikov conjecture originated from the work by Sergei Novikov on topological invariance of rational Pontryagin classes \cite{sN:64,sN:65a,sN:65b} that appeared in 1965.  The formal statement has been finalized by 1970, cf. \cite{sN:10,OB:95s}.  The history of this famous conjecture, which has its own Wikipedia page, where it is called ``one of the most important unsolved problems in topology'', is exposed in two volumes of proceedings \cite{OB:95} of the 1983 Oberwolfach conference and the book by Kreck and L\"{u}ck \cite{mKwL:05}.
We want to mention several important advances \cite{mBwHiM:93,mG:96,gY:98,aB:03,aD:04,eGnHsW:05,aDsFsW:08,sCsFgY:08,dRrTgY:11,gY:17} since that conference in proving the conjecture by different modern methods compared to the approach we explain here.
For an overview of the modern state of affairs we refer the reader to \cite{sW:90,jD:00,jR:16}. 

Continuing with the story, if the assembly is actually an injection then the signature is
a homotopy invariant.
This is the modern approach to proving the Novikov conjecture. 
In fact, stronger integral conjectures may be stated when integral
group homology is used, and there are K-,
A-theoretic, and $C^{\ast}$-algebraic analogues of these integral
maps.
For example, the statement parallel to ($\ddag$)  about classes in
$KO \lbrack {1\over2} \rbrack$ is equivalent to
integral injectivity of $\alpha$ \cite{sW:90}.
It makes sense, therefore, to call the conjecture
that the assembly map in K-theory is injective
for torsion-free $\Gamma$
{\it the integral Novikov conjecture in K-theory\/}.

Another interesting application can be obtained from
splitting the $C^{\ast}$-algebraic version of the assembly
map by applying the same approach
as taken here, giving what
J. Rosenberg calls the {\it strong Novikov conjecture\/}.
That is known to imply rigidity and vanishing results for
higher elliptic genera \cite{kL:92}.

In this paper we concentrate on the study of assembly maps in K-and L-theories relevant to the Borel conjecture on topological rigidity of closed aspherical manifolds. This connection is explained in detail in the next section. There are other versions of assembly maps in the literature and different ways to study them, most notably the Farrell-Jones assembly map defined by Davis and L\"{u}ck \cite{jDwL:98} and related conjectures about them.  We will discuss the relationships between these maps and conjectures in the last section \refS{fjcomp}.
For many other applications of the assembly maps we refer to the survey \cite{aBwLhR:08}.

\SSecRef{Topological rigidity of aspherical manifolds}{BRC}

We would like to emphasize one geometric application of the Loday assembly map. 

A prominent and long-standing conjecture of Armand Borel claims that any two homotopy equivalent closed aspherical manifolds (the manifolds with the contractible universal cover) are in fact homeomorphic.
This became known as the \textit{Topological Rigidity Conjecture} or \textit{Borel Conjecture}.
Similar relative rigidity results are stated for manifolds with boundary and even with corners.  The geometric precursor of this conjecture, which verifies the case of the Borel conjecture for hyperbolic manifolds, is Mostow's rigidity theorem.
For the history of this circle of ideas, we refer the interested reader to excellent survey chapters \cite{fF:02,cS:02}.

We include here a sketch of how the assembly maps can be used to prove this conjecture.

It is known that the Whitehead group $\textit{Wh} (\Gamma)$ contains obstructions to some crucial constructions used to classify manifolds $M$ with fundamental group $\Gamma$. We refer the reader to the following accessible literature.  The book by Rourke and Sanderson \cite{cRbS:72} for the definition of the Whitehead group in the context of piecewise-linear topology, where the crucial $s$-cobordism theorem is proved in chapter 6 as Theorem 6.19, or to a contemporary tutorial Hambleton \cite{iH:02}.  We mentioned the Whitehead group already in section \refSS{APPL} with relevant references to the literature and applications.  

Vanishing of the entire Whitehead group is a very much desired fact.  This group is the cokernel of the assembly map $A_{K,1} \colon H_1 (\Gamma, K(\mathbb{Z})) \to K_1 (\mathbb{Z} [\Gamma])$ induced from $\alpha_K$ in degree 1, see page 311 in Loday \cite{jlL:76}.  This means that $\textit{Wh} (\Gamma) = 0$ if $A_1$ is an epimorphism.  In this case, in particular, all $h$-cobordisms on $M$ are products.

A good way to rephrase this last fact is by using the notion of a structure set of $M$.  Consider the equivalence classes of maps of manifolds $f_{\ast} \colon M_{\ast} \to M$ under the equivalence relation $f_1 \sim f_2$ if there is a homeomorphism $h \colon M_1 \to M_2$ such that $f_2 h = f_1$.  This is the structure set $\mathcal{S} (M)$.
In these terms, the Borel Rigidity Conjecture is the claim that the cardinality of $\mathcal{S} (M)$ is 1 if $M$ is aspherical.
There is a variant of this notion, $\mathcal{S}^h (M)$, where the equivalence relation is weaker: instead of the existence of a homeomorphism $h$ one requires existence of an h-cobordism between $M_1$ and $M_2$.
Since nontrivial h-cobordisms over $M$ are in bijective correspondence with $\textit{Wh} (\Gamma)$, if $\textit{Wh} (\Gamma) =0$ then $\mathcal{S}^h (M) = \mathcal{S} (M)$.

Now the vanishing of $\mathcal{S}^h (M)$ can be addressed using the surgery long exact sequence.  For the background in algebraic L-theory the reader is encouraged to consult Ranicki \cite{aR:92}, where section 14 is devoted to the algebraic surgery exact sequence.  A portion of that sequence looks (for $M$ an aspherical manifold) as follows 
\begin{multline*}
\ \ \ \ \ \ \ H_{n+1} (M, L(\mathbb{Z})) \xrightarrow{\, A_{L,n+1}\, } L_{n+1} (\mathbb{Z} [\Gamma]) \longrightarrow \mathcal{S}^h (M)\\
\longrightarrow H_{n} (M, L(\mathbb{Z})) \xrightarrow{A_{L,n}} L_{n} (\mathbb{Z} [\Gamma]) \ \ \ \ \ \ \ 
\end{multline*}

If one is able to prove that the assembly maps $\alpha_K$ and $\alpha_L$ for the fundamental group of a specific aspherical manifold are weak equivalences, inducing isomorphisms in all dimensions, then one concludes that the Borel Conjecture is true in the case of that particular manifold.

In this paper we focus on the K-theoretic assembly map for groups $\Gamma$ with a finite model for $K(\Gamma, 1)$.
There are also likely extensions to groups that have torsion but with restrictions on the coefficient rings such as, for example, the rationals $\mathbb{Q}$. The Borel conjecture even with this constraint has interesting applications to other famous conjectures, for example Mathai's conjecture \cite{sC:04}.

There are analogues of our methods in L-theory, and the work on extending our results to L-theory is in progress.  

\SSecRef{Bounded K-theory and its application}{BBB}

We begin modeling both ends of the assembly map as fixed points of two $\Gamma$-equivariant spectra and the assembly itself as a fixed point map.
The interested reader can consult surveys of this material \cite{gC:93,gC:05} or the original papers in the references for details. 

\medskip

Let us assume we are given a finite complex which is a model for the classifying space $B\Gamma$.
The unversal cover $X$ of this complex is a cocompact model for the contractible free $\Gamma$-space $E\Gamma$.
The equivariant spectrum $S$ is then the locally finite homology of $X$ with coefficients in $KR$ denoted by $\hlf (X,KR)$.
It is a theorem \cite{gC:93} that $\hlf (X,KR)^{\Gamma} = B\Gamma_+ \wedge KR$.
Another theorem is that $r \colon \hlf (X,KR)^{\Gamma} \to \hlf (X,KR)^{h\Gamma}$ is an equivalence.

\medskip

For a proper metric space $X$ and a ring $R$,
Pedersen and Weibel \cite{ePcW:85}
define the category $\mathcal{B} (X,R)$ of \textit{geometric $R$-modules over $X$}.
The objects of this category are locally finite functions $F \colon X \to \Free_{fg} (R)$, from points of $X$ to the category of finitely generated free $R$-modules.  Following Pedersen and Weibel, we will denote by $F_x$ the module assigned to the point $x$ of $X$ and the object itself by writing down the collection $\{ F_x \}$.
The \textit{local finiteness} condition requires that for every bounded subset $S \subset X$ the restriction of $f$ to $S$ has finitely many nonzero modules as values.

Let $d$ be the distance function in $X$.  The morphisms $\phi \colon \{ F_x \} \to \{ G_x \}$ in $\mathcal{B} (X,R)$ are collections of $R$-linear homomorphisms
$\phi_{x,x'} \colon F_x \longrightarrow G_{x'}$,
for all $x$ and $x'$ in $X$, with the property
that $\phi_{x,x'}$ is the zero homomorphism whenever $d(x,x') > D$
for some fixed real number $D = D (\phi) \ge 0$.
In this case, we say that $\phi$ is \textit{bounded by} $D$.
The composition of two morphisms $\phi \colon \{ F_x \} \to \{ G_x \}$ and $\psi \colon \{ G_x \} \to \{ H_x \}$
is given by the formula
\[
(\psi \circ \phi)_{x,x'} = \sum_{z \in X} \psi_{z,x'} \circ \phi_{x,z}.
\]
The sum in the formula is finite because of the local finiteness of $\{ G_x \}$.

The category $\mathcal{B} (X,R)$ is certainly a symmetric monoidal category.
The associated K-theory is traditionally called the \textit{bounded K-theory}. 
It is clear that a free action on $X$ by isometries gives an induced free action on the category $\mathcal{B} (X,R)$.
What we need to do next is modify it into a useful equivariant theory.
For details we refer to \cite[ch.\,VI]{gC:95}.

\begin{NotRef}{OmitR}
$\mathcal{B} (X,R)$ can be thought of as a functor on the certain category of proper metric spaces $X$ and also a functor in the variable $R$.
For the purposes of this paper, $R$ can be fixed from the outset as an arbitrary ring in this and the next section.
The same is true about section \refS{SurjAss}, except $R$ will need to be noetherian.
In order to simplify the notation, we will use the notation $\mathcal{B} (X)$ for the category of geometric $R$-modules over $X$.
\end{NotRef}

Let $\EGamma$ be the category with the object set $\Gamma$ and the
unique morphism $\mu \colon \gamma_1 \to \gamma_2$ for any pair
$\gamma_1$, $\gamma_2 \in \Gamma$. There is a left $\Gamma$-action
on $\EGamma$ induced by the left multiplication in $\Gamma$.
If $\mathcal{C}$ is a small category with a left $\Gamma$-action,
then the category of functors
$\mathcal{C}_{\Gamma}=\Fun(\EGamma,\mathcal{C})$ is another
category with the $\Gamma$-action given on objects by
the formulas $\gamma(F)(\gamma')=\gamma F (\gamma^{-1} \gamma')$
and $\gamma(F)(\mu)=\gamma F (\gamma^{-1} \mu)$.
It is always
nonequivariantly equivalent to $\mathcal{C}$. The fixed
subcategory $\Fun(\EGamma,\mathcal{C})^{\Gamma} \subset
\mathcal{C}_{\Gamma}$ consists of equivariant functors and
equivariant natural transformations.

Explicitly, when $\mathcal{C} = \mathcal{B} (X)$ with the $\Gamma$-action described
above, the objects of $\mathcal{B}_{\Gamma} (X)^{\Gamma}$ are the pairs
$(F,\psi)$ where $F \in \mathcal{B} (X)$ and $\psi$ is a function on
$\Gamma$ with $\psi (\gamma) \in \Hom (F,\gamma F)$ such that
$\psi(1) = 1$, and $\psi (\gamma_1 \gamma_2) =
\gamma_1 \psi(\gamma_2)  \psi (\gamma_1)$.
These conditions imply that $\psi (\gamma)$ is always an
isomorphism as in \cite{rT:83}. The set of morphisms $(F,\psi) \to
(F',\psi')$ consists of the morphisms $\phi \colon F \to F'$ in
$\mathcal{B} (X)$ such that the squares
\[
\xymatrix{
F \ar[rr]^-{\psi (\gamma)} \ar[dd]_-{\phi} &&\gamma F \ar[dd]^-{\gamma \phi}\\
\\
F' \ar[rr]^-{\psi' (\gamma)} &&\gamma F'
}
\]
commute for all $\gamma \in \Gamma$. A slightly more refined
theory is obtained by replacing $\mathcal{B}_{\Gamma} (X)$ with the full
subcategory $\mathcal{B}_{\Gamma, 0} (X)$ of functors sending all
(iso)morphisms of $\EGamma$ such that the maps and their inverses are bounded by $0$. 
So $\mathcal{B}_{\Gamma, 0} (X)^{\Gamma}$ consists of $(F, \psi)$ with 
$\fil \psi (\gamma) = 0$ for all $\gamma \in \Gamma$.

The fixed point category 
$\mathcal{B}_{\Gamma, 0} (X)^{\Gamma}$ is clearly symmetric monoidal.
Now we have finally designed a spectrum 
$K (X)$ defined as the (nonconnective delooping of) the K-theory spectrum of
$\mathcal{B}_{\Gamma, 0} (X)$ 
with the desired property
\[
K (X)^{\Gamma} 
= K (\mathcal{B}_{\Gamma, 0} (X)^{\Gamma}) = K(R\Gamma).
\]

\medskip

We can now assemble the constructions and the facts into one commutative diagram as planned in the beginning of the section.
Here our choices are $S = \hlf (X,KR)$ and $T = K (X)$ so that $\underline{M} = K (X)^{h\Gamma}$.
\[
\xymatrix{
 B\Gamma_+ \wedge KR   \ar[dd]_-{\simeq} \ar[rr]^-{\alpha}
&&K(R\Gamma)  \ar[dd]^-{\simeq} \\
\\
 \hlf (X,KR)^\Gamma   \ar[dd]^-{r_S}_-{\simeq} \ar[rr]^-{\phi^{\Gamma}}
&&K (X)^\Gamma  \ar[dd]^-{r_T} \\
\\
 \hlf (X,KR)^{h\Gamma} \ar[rr]^-{\phi^{h\Gamma}}
&&K (X)^{h\Gamma}
}
\]

What is missing from the complete proof of the Novikov conjecture that $\alpha$ is a split injection is the definition of the equivariant \textit{controlled assembly map} 
$\phi \colon  \hlf (X,KR) \to K (X)$ and the verification that $\phi$ is a nonequivariant weak equivalence. 
Recall that this will automatically say that $\phi^{h\Gamma}$ is an equivalence.
The definition is straightforward; it is written down in \cite[Ch.\,4]{gC:93}.
Checking that $\phi$ is an equivalence is the whole point of the proof.
In fact, even before one deals with $\phi$, it would be very encouraging to realize similarities between 
$\hlf (X,KR)$ and $K (X)$.
The locally finite homology has classical origins, and the proper excision theorem is the main computational tool.
The bounded K-theory of Pedersen and Weibel is a much younger theory.
There is a more delicate \textit{bounded excision theorem} \cite[Ch.\,3]{gC:93} which can be matched to the excision in the locally finite homology in certain simple geometric situations (see the subsection below).
Luckily, excision computations in these simple geometric situations and the more elaborate derivations allow to realize that $\hlf (X,KR)$ and $K (X)$ are indeed equivalent.
It also follows that the controlled assembly $\phi$ is precisely a weak equivalence.

There are few examples in the literature where this primary strategy led to proofs of the integral Novikov conjecture.  The main difficulty is in applying the bounded excision in K-theory.  Still the size of the class of fundamental groups $\Gamma$ that can be handled this way is remarkable.  
The first application was to uniform lattices in connected semi-simple Lie groups in \cite{gC:93}. 
The authors proved in \cite{gCbG:04} that these groups have asymptotic dimension equal to the dimension of the symmetric space associated to the Lie group.  The proof of the Novikov conjecture was then extended to the class of all groups of finite asymptotic dimension \cite{gCbG:04}, the fact that was also known from the work of Bartels \cite{aB:03}.
These results are now subsumed by a more general theorem due to Ramras, Tessera, and Yu \cite{dRrTgY:11} which applies to the larger family of groups that have finite decomposition complexity.

\SSecRef{Example:\ free abelian groups}{BET}

The excision theorems are essential for computations needed in proving the Novikov conjecture, and their generalizations become the main theme in section \refS{SurjAss}.
Here we state them and apply to compute the assembly map for $\mathbb{Z}^n$.

Suppose $X$ is a proper metric space, in particular it is locally compact.
We also assume it is the union of countably many compact subsets.
Now let $X$ be an \textit{excisive} union of two closed subsets $U$ and $V$, so we assume that $X = \textrm{int} (U) \cup \textrm{int} (V)$.
There is the following expected excision theorem for the locally finite homology.

\begin{ThmRef}{Exc}
The commutative square
\[
\xymatrix{
 \hlf (U \cap V, KR) \ar[rr] \ar[dd]
&&\hlf (U,KR) \ar[dd] \\
\\
 \hlf (V,KR) \ar[rr]
&&\hlf (X,KR) }
\]
is a homotopy pushout.
\end{ThmRef}

\begin{proof}
The locally finite homology $\hlf$ here is called $\shlf$ in \cite[ch.\,II]{gC:93}.
This theorem follows from Proposition II.15 and Corollary II.17 in \cite{gC:93}.
\end{proof}

To state the correct analogue in bounded K-theory we will use the following notation.
Suppose $U$ is a subset of $X$.
Let $\mathcal{B}_{\Gamma, 0} (X)_{<U}$ denote the full subcategory of $\mathcal{B}_{\Gamma, 0} (X)$ on the objects $F$ with $F_x = 0$ for all points $x \in X$ with $d (x, U) \le D$ for some fixed number $D > 0$ specific to $F$.
This is an additive subcategory of $\mathcal{B}_{\Gamma, 0} (X)$ with the associated K-theory spectrum $K (X)_{<U}$.
It is easy to see that $\mathcal{B}_{\Gamma, 0} (X)_{<U}$ is in fact equivalent to $\mathcal{B}_{\Gamma, 0} (U)$, so $K (X)_{<U}$ is equivalent to $K (U)$.
Similarly, if $U$ and $V$ are a pair of subsets of $X$, then there is the full additive subcategory
$\mathcal{B}_{\Gamma, 0} (X)_{<U,V}$ of $F$ with $F_x = 0$ for all $x$ with $d (x, U) \le D_1$ and $d (x, V) \le D_2$ for some numbers $D_1, D_2 > 0$.
It is an important observation that $\mathcal{B}_{\Gamma, 0} (X)_{<U,V}$ is not in general equivalent to $\mathcal{B}_{\Gamma, 0} (U \cap V)$.

\begin{ThmRefName}{BddExc}{Bounded Excision}
For a proper metric space $X$ and a pair of subsets $U$, $V$ which form a cover of $X$, the commutative square
\[
\xymatrix{
 K (X)_{<U, V} \ar[rr] \ar[dd]
&&K (U) \ar[dd] \\
\\
 K (V) \ar[rr]
&&K (X) }
\]
is a homotopy pushout.
\end{ThmRefName}

This is a consequence of Corollary IV.6 in \cite{gC:93}.

Finally, there is a situation when $K (X)_{<U, V}$ and $K (U \cap V)$ are equivalent.

\begin{DefRef}{CATHJK}
A pair of subsets $U$, $V$ of $X$ is called \textit{coarsely antithetic} if for each number $K > 0$ there is a number $K' > 0$ so that $U[K] \cap V[K] \subset (U \cap V)[K']$.
\end{DefRef}

Examples of coarsely antithetic pairs
include any two non-vacuously intersecting closed subsets of a simplicial tree as well as
complementary closed half-spaces in a Euclidean space.

\begin{CorRef}{BDDEXCII}
If $U$ and $V$ is a coarsely antithetic pair of subsets of $X$ which form a cover of $X$, then the commutative square
\[
\xymatrix{
 K (U \cap V) \ar[rr] \ar[dd]
&&K (U) \ar[dd] \\
\\
 K (V) \ar[rr]
&&K (X)
}
\]
is a homotopy pushout.
\end{CorRef}

We will use these theorems to compute the map 
$\phi \colon \hlf (\mathbb{Z}^n,KR) \to K(E\mathbb{Z}^n)$.
The $n$-dimensional torus $T^n$ is a model for $B \mathbb{Z}^n$.
Then $E \mathbb{Z}^n$ is $\mathbb{R}^n$.
Consider the covering of $\mathbb{R}^n$ by two subsets, the upper half-space $\mathbb{R}^n_+$ and the lower half-space $\mathbb{R}^n_-$.  Applying Theorem \refT{Exc} to this covering and using the fact that both covering sets are $\hlf$-acyclic, we see from the Mayer-Vietoris sequence that 
$\hlf (\mathbb{R}^n,KR) \simeq \Sigma \hlf (\mathbb{R}^{n-1},KR)$.
Inductively, we obtain $\hlf (\mathbb{R}^n,KR) \simeq \Sigma^n \hlf (\mathbb{R}^{0},KR)$, where $\mathbb{R}^{0}$ is the point at the origin of $\mathbb{R}^n$, so $\hlf (\mathbb{R}^{0},KR) = KR$.

Observing that $\mathbb{R}^n_+$ and $\mathbb{R}^n_-$ form a coarsely antithetic covering, we can apply Corollary \refC{BDDEXCII} in the same fashion to compute 
$K (\mathbb{R}^n) \simeq \Sigma^n K (\mathbb{R}^{0}) \simeq \Sigma^n KR$.

To see that $\phi$ in fact induces the weak equivalence we detected, we need to examine the assembly maps on the terms of the homotopy colimits.
The components are the assembly maps $\phi^k_{\pm} \colon \hlf (\mathbb{R}^k_{\pm},KR) \to K(\mathbb{R}^k_{\pm})$, which map a contractible spectrum to a contractible spectrum, and 
$\phi^k \colon \hlf (\mathbb{R}^k,KR) \to K(\mathbb{R}^k)$ for $0 \le k < n$, which we assume are equivalences by induction.  This gives us the equivalence of the homotopy pushouts, which is $\phi$.

\SecRef{Applications of continuous control at infinity}{reen}

\SSecRef{Continuous control in good compactifications}{nice}

Continuous control at infinity was introduced in \cite{gCeP:93}.
We will take a narrower point of view and assume that we are given a metric space $X$.
Now suppose $X$ is an open dense subset of a compact Hausdorff space $\widehat{X}$, and the difference is $Y = \widehat{X} \setminus X$.
The objects $F$ of the category $\mathcal{C} (\widehat{X},Y)$ are the same as those of $\mathcal{B} (X)$.
If $U \subset \widehat{X}$ is a subset,
and $F \in {\cal B} (X)$,
define $F \vert U$ by $(F \vert U)_x = F_x$
if $x \in U \setminus Y$ and 0 if $x \in X \setminus U \setminus Y$.
An $R$-homomorphism $\phi \colon F \rightarrow G$ in
is called {\it continuously controlled\/} at $y \in Y$
if for every neighborhood $U$ of $y$ in $\widehat{X}$
there is a neighborhood $V$ so that
$\phi (F \vert V) \subset G \vert U$ and
$\phi (F \vert (\widehat{X} \setminus U)) \subset G \vert (\widehat{X} \setminus V)$.
The morphisms of $\mathcal{C} (\widehat{X},Y)$ are the homomorphisms which are continuously
controlled at all $y \in Y$.

Now let $W$ be an open subset of $\widehat{X}$ and
$p\colon W \rightarrow K$ be a map with continuous restriction
$p \vert (Y \cap W)$.
A morphism $\phi\colon F \rightarrow G$ in $\mathcal{C} (\widehat{X},Y)$
is {\it $p$-controlled\/} at $y \in Y \cap W$
if for every neighborhood $U$ of $p(y)$ in K
there is a neighborhood $V$ of $y$ in $\widehat{X}$
so that $\phi (F \vert V) \subset G \vert p^{-1}(U)$
and $\phi \bigl( F \vert (\widehat{X} \setminus p^{-1}U) \bigr) \subset B \vert (\widehat{X} \setminus V)$. 
The category $\mathcal{C} (\widehat{X},Y,p)$ has the same objects as
$\mathcal{B} (X)$ and morphisms which are $R$-linear homomorphisms continuously controlled
at all $y \in Y \setminus W$ and $p$-controlled at all $y \in W \cap Y$.

Both $\mathcal{C} (\widehat{X},Y)$ and $\mathcal{C} (\widehat{X},Y,p)$  are small symmetric monoidal categories.
We still assume there is a free action of $\Gamma$ on $X$ by isometries such that the orbits space is a finite model for $B\Gamma$.
If this action of $\Gamma$ on $X$ extends continuously to $Y$, and if the extension restricts to an action on $W$, then there are induced actions.
In fact, there are better equivariant theories as before.
We will use the notation $K (\widehat{X},Y)$ and $K (\widehat{X},Y,p)$ for the resulting equivariant spectra.

We can still use the new models to interpret the same assembly map $\alpha$.

\medskip

Let $C \widehat{X}$ be the cone on $\widehat{X}$ with
the identification $\widehat{X} = \widehat{X} \times \{ 1 \}$
and let $p\colon  \widehat{X} \times (0,1) \rightarrow \widehat{X}$ be the projection.
The map $\pi\colon C \widehat{X} \rightarrow \Sigma \widehat{X}$ collapsing the base 
$\widehat{X}$ induces an equivariant functor
\[
\pi_\ast \colon \mathcal{C}(C \widehat{X},CY \cup \widehat{X},p)
\to
\mathcal{C}(\Sigma \widehat{X}, \Sigma Y,p)
\]
which in turn induces a map of spectra
\[
\pi_\ast \colon K(C \widehat{X},CY \cup \widehat{X},p)
\longrightarrow
K (\Sigma \widehat{X}, \Sigma Y,p).
\]

In the remainder of this section we assume that $X$ is the universal cover of a finite complex K with $\pi_1 (K) = \Gamma$.
If we choose $S = \Omega K(C \widehat{X},CY \cup \widehat{X},p)$ and $T = \Omega K (\Sigma \widehat{X}, \Sigma Y,p)$ then 
there is a commutative diagram
\[
\xymatrix{
B\Gamma_+\wedge KR \ar[rr]^-{\alpha} \ar[dd]_-{\simeq} &&K(R\Gamma) \ar[dd]_-{\simeq}   \\
\\
S^{\Gamma}  \ar[rr]^-{\Omega \pi_{\ast}^{\Gamma} }  &&T^{\Gamma}
}
\]
The vertical equivalences and the fact that the canonical map 
$r_S \colon {S}^{\Gamma} \to {S}^{h\Gamma}$ is an equivalence
are verified in \cite{gCeP:93}.
Again, in order to prove instances of Novikov conjecture by choosing the target 
$\underline{M} = {T}^{h\Gamma}$, one needs to look at specific geometric situations where the equivariant map $\phi \colon S \to T$ is an nonequivariant equivalence.
The following is a list of general conditions from \cite{gCeP:93} that accomplish this:
\begin{enumerate}
\renewcommand{\labelenumi}{(\roman{enumi})}
\item The $\Gamma$-action
extends continuously to $\widehat{X}$,
\item $\widehat{X}$ is metrizable,
\item $\widehat{X}$ is contractible,
\item compact subsets of $X$ become small near $Y$ or, more precisely,
for every point $y$ in $Y$, for every compact subset $K \subset X$ and
for every  neighborhood $U$ of $y$ in $\widehat{X}$, there exists a neighborhood $V$ of
$y$ in $\widehat{X}$ so that if $g\in \Gamma$ and $gK\cap V \ne \emptyset$ then $gK
\subset U$.
\end{enumerate}
We will refer to such compactifications of $X$ as \textit{good} when they are available.
The last property (iv) is usually expressed by saying that the action of $\Gamma$ on $X$ is \textit{small at infinity}.

This kind of situation is remarkably common in geometry.
There are many classical examples of universal covers $X$ of compact aspherical manifolds that satisfy this list of conditions.
The most common and very general curvature condition that gives these properties is CAT(0).
It also takes us outside of the world of spaces of finite asymptotic dimension that we worked out in section 2.

\SSecRef{Continuous control with large actions at infinity}{PfTh2}

We move away from the development of general conditions and theorems and review more ad hoc methods of constructing proofs of Novikov conjecture.

Given that $r_S$ is an equivalence, we should try to manufacture the homotopy type of $S$ in the form of an equivariant spectrum $M$.
The point is that this time we will not insist on the homotopy fixed point being equivalent to $T^{h\Gamma}$.
What we do need is that $M$ receives an equivariant map $\beta$ from $T$ and the homotopy equivalence we are designing is realized by the composition $\beta \circ \phi$.
Then the following diagram gives us a splitting of $f$ according to the same principles
\[
\xymatrix{
 S^\Gamma   \ar[dd]^-{r_S} \ar[rr]^-{f}
&&T^\Gamma  \ar[dd]^-{r_T} \ar@/^1pc/[dddr] \\ 
\\
 S^{h\Gamma} \ar[rr]^-{\phi^{h\Gamma}} \ar@/_1pc/[rrrd]
&&T^{h\Gamma} \ar[dr]^{\beta^{h\Gamma}}\\
&&&M^{h\Gamma}
}
\]
The new target $\underline{M}$ is simply further removed from the assembly map.
The advantage of the new target is its simplicity and cusomizability.

Consider the map $\kappa$ with domain $C \widehat X$
which contracts the subspace $CY$ and
produces the reduced cone $\widetilde C X^+$.
It induces an equivariant functor
\[
\lambda \colon
\mathcal{C}(C \widehat X ,CY \cup \widehat X,p)
\longrightarrow
\mathcal{C}(\widetilde C X^{+}, X^+)
\]
because each morphism from
$\mathcal{C}(C \widehat X , CY \cup \widehat X,p)$
is controlled at $X^{+}$.
(In other words, the $p$-control is automatic.)
This functor gives an equivariant map of spectra
\[
\lambda \colon
{S}=\Omega K(C \widehat X ,CY \cup \widehat X,p)
\longrightarrow
{Q} \overset{ \text{def} }{=} \Omega K(\widetilde C X^+, X^+).
\]
which is a weak homotopy equivalence \cite{gCeP:98}.

Since $X^+$ is metrizable, ${Q}$ is a Steenrod functor by
\cite[Theorem 1.36]{gCeP:93}.
Also ${Q}^{\Gamma} \simeq {Q}^{h\Gamma}$ as before.
Another Steenrod functor is the \v{C}ech homology $\check h (X^+;KR)$ defined in \cite{gCeP:98, bG:97} as follows.

\begin{DefRef}{CechHom}
The {\it \v{C}ech homology\/} of a space $Z$ with
coefficients in $S$ is the (simplicial)
spectrum valued functor
\[
\check h (Z;S) = \holim{Cov \, Z} N \underline {\phantom {S}} \wedge S,
\]
where $Cov \, Z$ is the category of finite rigid open coverings
of $Z$ defined in \cite{gCeP:98} and $N$ is the (simplicial) nerve.
This is a generalized Steenrod homology theory.
\end{DefRef}

The {\it support at infinity\/} of an object
$F \in \mathcal{C}(\widehat X,Y)$ is the set of limit points of
$\{ x \vert F_x \ne 0 \}$ in $Y$.
The full subcategory of $\mathcal{C}(\widehat X,Y)$
on objects with support at infinity contained in
$C \subset Y$ is denoted by $\mathcal{C}(\widehat X,Y)_{C}$.

The inclusions of two open sets $U_1, U_2$ in $X^+$
induce the maps
\[
K(\widetilde C X^+, X^+)_{U_1 \cap U_2}
\longrightarrow 
K(\widetilde C X^+, X^+)_{U_i}.
\]
In general, there is a functor
$Int \, \mathcal{U} \rightarrow spectra$ for any
$\mathcal{U} \in Cov \, X^+$,
where $Int \, \mathcal{U}$
is the partially ordered set of
all multiple finite intersections
of members of $\mathcal{U}$.

There is the following controlled excision result.

\begin{ThmRef}{Claimtwo}
For a fixed $\mathcal{U} \in Cov \, X^+$,
the universal excision map
\[
\hocolim{Int \, \mathcal{U}} K(\widetilde C X^+, X^+)_{\cap U_i}
\longrightarrow
K(\widetilde C X^+, X^+)
\]
is a weak equivalence.
\end{ThmRef}

The spectrum $\Sigma {Q}$ on the right is a $\Gamma$-equivariant spectrum.
To rediscover this aspect of the structure on the left,
we can write
\[
\holim{\mathcal{U} \in Cov X^+}
\left( \hocolim{Int \, \mathcal{U}} K(\widetilde C X^+, X^+)_{\cap U_i} \right)
\simeq \Sigma {Q},
\]
where the $\Gamma$-action on the left-hand side
is induced from the obvious action on $Cov \, X^+$.
By sending each non-empty $\cap U_i$ to a point,
one gets maps
\[
\hocolim{Int \, \mathcal{U}} K(\widetilde C X^+, X^+)_{\cap U_i}
\longrightarrow
\vert Int \, \mathcal{U} \vert \wedge KR.
\]
Finally, we get the induced equivariant map
of homotopy limits
\[
\pi \colon  {Q} \longrightarrow \check h (X^+;KR).
\]
This map can be viewed as a component of a natural transformation
of Steenrod functors which is an equivalence on the point,
hence $\pi$ is an equivalence,
according to a theorem of Milnor.
Recall that this is enough to see that
${Q}^{h\Gamma} \simeq \check h (E\Gamma^+;KR)^{h\Gamma}$.

Recall the notation $S = \Omega K(C \widehat{X},CY \cup \widehat{X},p)$ and $T = \Omega K (\Sigma \widehat{X}, \Sigma Y,p)$.
The following diagram summarizes what has been accomplished so far, before we start making choices for $\underline{M}$.
\begin{equation}
\begin{gathered}
\xymatrix{
 S^\Gamma   \ar[dd]^-{r_S}_-{\simeq} \ar[rr]^-{\alpha}
&&T^\Gamma  \ar[dd]^-{r_T} \\
\\
 S^{h\Gamma} \ar[rr]^-{\phi^{h\Gamma}} \ar[dd]^-{\pi^{h\Gamma} \circ \lambda}_-{\simeq}
&&T^{h\Gamma} \ar[dd] \\
\\
 \check h (X^+;KR)^{h\Gamma}  \ar[rr]
&&\underline{M}  
}
\end{gathered}
\tag{$\dagger$}
\end{equation}

\medskip

Returning to ${\cal T}$, we would like an excision result
analogous to Theorem \refT{Claimtwo}.
One important point is that there is no automatic $p$-control as in
the construction of $\lambda$ because
the action is no longer assumed to be small at infinity.
This may not allow the excision map to exist at all.
In order to produce such a map, the covering sets $p(U) \subset Y$
need to be boundedly saturated.

\begin{DefRef}{BSDef}
For any subset $K$ of a metric space $(X,d)$
let $K[D]$ denote the set $\{ x \in X : d(x,K) \le D \}$ called the $D$-\textit{neighborhood} of $K$.

If a metric space $(X,d)$ is embedded in a topological
space $\widehat X$ as an open dense subset,
a set $A \subset Y =\widehat X \setminus X$ is
{\it boundedly saturated\/} if for every closed subset
$C$ of $\widehat X$ with $C \cap Y \subset A$, the closure
of each $D$-neighborhood of $C \backslash Y$ for
$D \ge 0$ satisfies
$\overline{(C \backslash Y)[D]} \cap Y \subset A$.
\end{DefRef}

We obtain a map
\[
\pi_{T} \colon  {\cal T}
\longrightarrow
\Sigma \holim{\alpha \in \{ \alpha \}}
N\alpha \wedge KR
\]
for any $\Gamma$-invariant system of coverings of $Y$
by boundedly saturated open sets.
This map is $\Gamma$-equivariant,
so the composition induces a map
\[
{\cal T}^{h\Gamma} \longrightarrow
\left( \Sigma \holim{\alpha \in \{ \alpha \}}
N\alpha \wedge KR \right)^{\!\!h\Gamma} \overset{ \text{def} }{=} \underline{M}.
\]
This map may exist even if the boundedly saturated sets in $\{ \alpha \}$ are not necessarily open.
If there is a $\Gamma$-invariant category $\{ \alpha \}$ of finite coverings $\alpha$ by boundedly saturated sets in $Y$ so that there is an analogue of the excision result from \cite{gCeP:98}, 
we call $\{ \alpha \}$ \textit{excisive}.

How does one generate a good $\alpha$?
It is clear that the boundedly saturated sets are closed under union.
One option is to ``saturate'' the open sets $U \subset Y$
by associating to $U$ its envelope in 
a Boolean algebra of all
boundedly saturated subsets of $Y$
thus mapping $Cov \, Y$ functorially onto
the resulting category $\{ \alpha \}$.
Let us denote this functor by 
$\textit{sat} \colon \beta \mapsto \alpha (\beta)$.
Collections of open saturated sets are excisive, but 
the saturation process does not preserve openness.
However, one can often make an
intelligent choice of the Boolean subalgebra of
boundedly saturated sets in the construction of $\{ \alpha \}$.
In authors' experience, once the relevant specific features of the asymptotic geometry of $X$ are understood,
such a choice for the excisive subalgebra $\{ \alpha \}$ becomes natural.
Let us assume for simplicity that $\{ \alpha \}$ are open coverings.

With the established target $\underline{M}$, the next task is to connect
$\check h (X^+;KR)^{h\Gamma} \to \underline{M}$.
There is a map
\[
\theta\colon \check h (Y;KR)
\longrightarrow \holim{\{ \alpha \}} N
\underline {\phantom {S}} \wedge KR
\]
induced by the inclusion of categories
$\{ \alpha \} \hookrightarrow Cov \, Y$.
Assuming $\widehat X$ is \v Cech-acyclic,
there is a composite weak equivalence
\[
\xymatrix{
\holim{Cov \widehat X \cup CY} \!\!\! N
\underline {\phantom {S}} \wedge KR 
\ar[rr]^-{\simeq} \ar[dd]^-{\simeq} 
&&\holim{Cov \, \Sigma Y} N
\underline {\phantom {S}} \wedge KR \ar[dd]_-{\simeq} \\
\\
\holim{Cov X^+} N
\underline {\phantom {S}} \wedge KR  
\ar[rr]^-{\exists}  &&\Sigma \holim{Cov \, Y} N
\underline {\phantom {S}} \wedge KR
}
\]
Now the bottom of diagram ($\dagger$) factors as 
\[
\xymatrix{
 \!\!\!\!\!\!\!\!\!\!\!\!\!\!\!\!\!\!\!\! 
 \check h (X^+;KR)^{h\Gamma}  \ar[rr] \ar@/_1pc/[ddr]_-{\simeq}
&&\underline{M}  \\
\\
&\Sigma \check h (Y;KR)^{h\Gamma} \ar@/_1pc/[uur]_-{\Sigma \theta^{h\Gamma}}
}
\]

As before, if $\theta$ is a weak equivalence then
$\Sigma \theta^{h\Gamma}$ is a weak equivalence, so
$\alpha$ is a split injection.
Of course, this is a major "if" but it is also an opportunity to select 
$\{ \alpha \}$ to get the homotopy type of 
$\holimplain_{\{ \alpha \}} N
\underline {\phantom {S}} \wedge KR$ just right:
it has to be equivalent to the domain of $\theta$ which is the \v{C}ech homology spectrum $\check h (Y;KR)$.

\begin{Sum}
Given a discrete group $\Gamma$, the method
described here calls for a construction of a
compact classifying space $B \Gamma$ and
an equivariant compactification $\widehat X$ of the universal
cover $E \Gamma$.
The space $\widehat X$ itself may not be metrizable but it
is required to be acyclic in the sense that
its \v Cech homology is that of a point.
Then a convenient metric must be introduced on $E \Gamma$.
The action on $\widehat X$ may not be small at infinity,
but the choice of metric determines the family of
boundedly saturated subsets of the boundary $Y=\widehat X \setminus E \Gamma$.
One has to make a choice of a $\Gamma$-invariant
collection of open coverings of $Y$ by such sets
which preserves the \v Cech homotopy type of $Y$.
Furthermore, the weak homotopy equivalence of
\v Cech homology spectra has to be realized by
the map $\theta$ defined above.
\end{Sum}

\begin{Cav}
This is surely a generaliation of the strategy in section \refSS{nice}.
If the compactification is small at infinity, all sets are boundedly saturated, so the easy choice of $\{ \alpha \}$ would be all all finite open coverings making $\theta$ the identity map.
In more general situations, the collection of open boundedly saturated subsets in $Y$ is usually too coarse to reproduce the \v{C}ech homotopy type of the boundary.
The better choices use excisive boundedly saturated coverings which are not open.

There is often a convenient cofinal
category $C$ of $Cov \, Y$
where the open covering sets have predictable saturation in $Y$
which (1) does not change the homotopy type of the nerve for each individual covering,
in other words the inclusion $N \beta \hookrightarrow N \alpha (\beta)$ is a homotopy equivalence for each $\beta$, and (2) produces boundedly saturated coverings $\{ \alpha' \}$ that are still excisive.

In this case, the map $\theta$ can be reconstucted as follows.
Since $\textit{sat}$ is left cofinal,
and the assignment $\beta \mapsto \alpha (\beta)$
induces a natural transformation
of the functors
$N \beta \wedge KR \rightarrow N \alpha (\beta) \wedge KR$
from $Cov \, Y$ to spectra,
we can induce the following maps

\bigskip

\[
\begin{CD}
\holim{\beta \in Cov \, Y} \!\! N \beta \wedge KR
@>{\mathit{sat}_{\ast}}>>
\holim{\beta \in Cov \, Y} \!\! N \alpha (\beta) \wedge KR
@<{\simeq}<<
\holim{\alpha \in \{ \alpha \}} N \alpha \wedge KR \\
@V{\imath^{\ast}}VV @VVV @VVV \\
\holim{{C}} N \underline{\phantom{S}} \wedge KR
@>{(\mathit{sat} \vert {C})_{\ast}}>>
\holim{{C}} N \alpha (\underline{\phantom{S}}) \wedge KR
@<{\simeq}<<
\holim{\{ \alpha' \}} N \underline{\phantom{S}} \wedge KR.
\end{CD}
\]

\bigskip

Here the vertical maps are induced by inclusions.
Now the map
${T} \rightarrow \holimplain N \alpha \wedge KR$
can be composed with the vertical map on the right,
so in order to split the assembly map we need
$\imath^{\ast}$ and $(\mathit{sat} \vert {C})_{\ast}$
to be weak equivalences.
This is ensured by the cofinality of $C$ and its special saturation properties.
\end{Cav}

\begin{ExRefName}{CAT0}{CAT(0) groups}
Let $X$ be a complete metric space satisfying the CAT(0) property.
The CAT(0) spaces generalize classical geometric objects such as symmetric spaces and buildings associated to classical semi-simple Lie groups.
The divergency properties of the geodesics in $X$ can be used to show that these spaces are always contractible, so $X$ is indeed a model for $E\Gamma$ for any cocompact lattice acting freely on $X$.
Another consequence of the CAT(0) condition is the fact that equivalence classes of geodesic rays can be organized into the \textit{ideal boundary} $\partial X$ so that the resulting compactification $\widehat{X} = X \cup \partial X$ is small at infinity, is metrizable, and is contractible.
In other words, this is a good compactification, so the assembly map $\alpha$ for such CAT(0)-groups $\Gamma$ is split injective by section \refSS{nice}.
\end{ExRefName}

\begin{ExRefName}{RHL}{Relatively hyperbolic lattices}
Let $X$ is the negatively curved symmetric space
associated to a semi-simple linear algebraic $\mathbb{Q}$-group
of real rank one and $\Gamma$ is a torsion-free arithmetic subgroup.
More generally, one can consider a complete non-compact
finite-volume Riemannian manifold $M$ with pinched negative
sectional curvatures $-a^2 \le K \le - b^2 <0$ and
torsion-free fundamental group $\Gamma$.
The universal cover $X$ of $M$ is hyperbolic relative to
the cusp subgroups.

For simplicity, assume that there is a unique cusp of $V$ with
the corresponding cusp subgroup $H$.
It is known that $H$ is a torsion-free finitely generated
nilpotent group.
The stratum, or horosphere, $X_H$ can be viewed as the underlying space
of a connected simply connected nilpotent Lie group where
$H$ acts by left multiplication.

We start with a construction of $\widehat{X}_H$ which is quite simple.
Geometrically $X_H$ is a Euclidean space but the compactification of 
$X_H$ by the ideal boundary is not $H$-equivariant.
Instead, one can use the ``box'' compactification with some collapses in the faces.  The details can be seen in \cite{bG:97}.
In fact, one can obtain a bonus proof of the Novikov conjecture for $H$
using this construction.  Recall that in particular this means we have a collection of boundedly saturated sets that give coverings $C_H$.

The points of $X_H$ are represented by geodesics converging to the same ideal point in $\partial X$.
The isometric action by $\Gamma$ allows to identify other ideal points, which we 
now resolve into copies of $X_H$.
Each of these copies gets compactified as above.  
The boundary $Y$ is the union of the resolutions and the rest of the ideal points.
There is a compact topology on $\widehat{X}$ which has all the required properties for the proof of Novikov.
In particular, it turns out that each compactified stratum and each of the remaining ideal points are boundedly saturated subsets of $Y$
The coverings like $C_H$ in each stratum and the ``geodesic shadows'' they generate on $\partial X$ give a family of boundedly saturated sets in $Y$ and form the needed family of coverings $C$.  We refer for details to \cite{bG:97}.
\end{ExRefName}

\begin{RemRef}{YICD}
This argument has been refined to apply to higher rank arithmetic groups \cite{bG:99} and to general relatively hyperbolic groups in the sense of Gromov \cite{bG:98}.
It can also be used in much fancier situations where the geometry of $X$ is understood only to a limited extent.
For example, it can be applied to the mapping class groups and their moduli spaces.
The mapping class groups have finite asymptotic dimension \cite{mBkBkF:10},
so this case of Novikov can also be proven by the methods in section \refS{bddkth}.
However this last example should underscore the flexibility and ad hoc nature of the method.
It is very likely that we already know enough about the Outer space in order to prove Novikov for groups of automorphisms of free groups.
\end{RemRef}

\SecRef{Surjectivity of the K-theoretic assembly map}{SurjAss}

The strategy for proving surjectivity of the Loday assembly map $\alpha_K$ that was suggested in the introduction is to split the splitting of $\alpha$ in those cases when that has been shown to exist.
In this section we assume that $R$ has finite global dimension (see the review of related algebra in section \refSS{CCFG}).
We will also assume that the group $\Gamma$ is \textit{geometrically finite}, 
in the sense that it has finite $B\Gamma$.

\begin{ThmRefName}{MTBB}{Main Theorem of \cite{gCbG:13}}
Suppose $\Gamma$ is a geometrically finite group of finite asymptotic dimension and $R$ is a regular noetherian ring of finite global dimension.
Then the assembly map
\[
\alpha_K \colon
B\Gamma_+\wedge K (R) \longrightarrow K (R\Gamma)
\]
is an equivalence.
\end{ThmRefName}

In this case, there is a splitting $r_T \colon K(X)^{\Gamma} \to K(X)^{h\Gamma}$ constructed in section \refS{bddkth}, where $X$ is a model for $E\Gamma$.
The goal of this section is to sketch the proof that $r_T$ is split injective.
The same strategy as before asks for a target $\underline{M}$ to fit a commutative diagram
\[
\xymatrix{
 K(X)^{\Gamma}   \ar[ddr]^-{s} \ar[rr]^-{r_T}
&&K(X)^{h\Gamma}  \ar[ddl]_-{t} \\
\\
& \underline{M}
}
\]
so that $s$ is an equivalence.
We will be working on a sequence of constructions starting with $K(X)^{h\Gamma}$ that allow to further extend the map $t$ while moving toward the final target $\underline{M}$ that we are able to guarantee is equivalent to $K(X)^{\Gamma}$.
This time however $\underline{M}$ will look like and smell like but not be a fixed point spectrum.

\SSecRef{Preparation for the argument}{gjudfj}

We can assume that $\Gamma$ is the fundamental group of a compact
aspherical manifold $M$ of dimension $n$, possibly with boundary, and that $X$ is the universal cover of $M$.
There is a proper metric on $X$ which is commensurable with the given word metric on $\Gamma$ when viewed as a free orbit in $X$.
Suppose the manifold $M$ is embedded in a sphere
$S^{n+k}$ for sufficiently large $k$. Let $Y$ be the universal
cover of the small tubular neighborhood $N$ of the embedding.
There is a metric on $Y$ commensurable with the metric on $X$.
Let $\partial Y$ denote the topological boundary of $Y$.

There are two geometric constructions and one variant of bounded K-theory in \cite{gCbG:13} that we need to review.

\begin{DefRef}{TX}
Starting with any set $Z$, let $S \subset Z \times Z$ denote any symmetric and
reflexive subset with the property that
\begin{itemize}
\item for any $z$, $z'$, there are elements $z_0$, $z_1$, ... , $z_n$ so that $z_0 = z$, $z_n = z_0$,
and $(z_i, z_{i+1}) \in S$.
\end{itemize}
Let $\rho \colon S \to \mathbb{R}$ be any function for which the following
properties hold:
\begin{itemize}
\item $\rho (z_1, z_2) = \rho (z_2, z_1)$ for all $(z_1, z_2) \in S$,
\item $\rho (z_1, z_2) = 0$ if and only if $z_1 = z_2$.
\end{itemize}
Given such $S$ and $\rho$, we define a metric $d$ on $Z$ to be the largest metric $D$
for which $D (z_1, z_2) \le \rho (z_1, z_2)$ for all $(z_1, z_2) \in S$. This means that $d$ is given
by
\[
d(z_1, z_2) = \inf_{n,\{z_0,z_1,...,z_n\}}
\sum_{i=0}^n
\rho (z_i, z_{i+1}).
\]

For any metric space $X$, the new metric space $TX$ is related to the cone
construction.
The underlying set is $X \times \mathbb{R}$. 
Let $S$ to be the set consisting of pairs of
the form $((x, r), (x, r'))$ or of the form $((x, r), (x', r))$. 
Then define $\rho$ on $S$
by
\[
\rho ((x, r), (x, r')) = \vert r - r' \vert ,
\]
and
\[
\rho ((x, r), (x', r)) = \left\{
                           \begin{array}{ll}
                             d (x, x'), & \hbox{if $r \le 1$;} \\
                             r d (x, x'), & \hbox{if $1 \le r$.}
                           \end{array}
                         \right.
\]

Since $\rho$ clearly satisfies the hypotheses of the above definition, we set the metric
on $T X$ to be $d$.
We also extend the definition to pairs of metric spaces $(X, Y)$,
where $Y$ is given the restriction of the metric on $X$.
\end{DefRef}

It is time to discuss functoriality in bounded K-theory.
Let $X$ and $Y$ be proper metric spaces with distance functions $d_X$ and $d_Y$.
A map $f \colon X \to Y$
of proper metric spaces is \textit{eventually continuous} if there is a real positive function $l$
such that
\begin{equation}
d_X (x_1, x_2) \le r \ \Longrightarrow
d_Y (f(x_1), f(x_2)) \le l(r).  \tag{1}
\end{equation}

A map $f \colon X \to Y$
of proper metric spaces is \textit{proper} if $f^{-1} (S)$ is a bounded subset of $X$ for
each bounded subset $S$ of $Y$.
We say $f$ is a \textit{coarse map} if it is proper and eventually continuous.
For example, all \textit{bounded} functions $f \colon X \to X$ with
$d_X (x, f(x)) \le D$, for all $x \in X$ and a fixed $D \ge 0$, are coarse.

The map $f$ is a \textit{coarse equivalence} if there is a coarse map $g \colon Y \to X$ such
that $f \circ g$ and $g \circ f$ are bounded maps.

Coarse maps are precisely the maps between proper metric spaces that induce maps on bounded K-theory spectra.

Another feature of bounded K-theory is that its objects can be generalized from $R$-modules locally modeled on finitely generated free $R$-modules to objects locally modeled on an arbitary small additive category $\mathcal{A}$.
One simply makes stalks $F_x$ be objects from the category $\mathcal{A}$.
We will use notation $\mathcal{B}(X, \mathcal{A})$ for such a theory.

\begin{DefRefName}{POI1}{Coarse Equivariant Theories}
We associate two new equivariant theories to metric spaces with a $\Gamma$-action.  The theory $K_i^{\Gamma}$ is defined only for metric
spaces with actions by isometries, while $K_p^{\Gamma}$ only for metric
spaces with coarse actions.

\smallskip

$K^{\Gamma}_i (Y)$ is the nonconnective $K$-theory of
$\mathcal{B}^{\Gamma}_i (Y) = \mathcal{B}_{\Gamma,0}(\Gamma \times Y)$,
where $\Gamma$ is regarded as a
word-length metric space with isometric $\Gamma$-action given by left multiplication, and $\Gamma
\times Y$ is given the product metric and the product isometric action.

\smallskip

$K^{\Gamma}_p (Y)$ is defined for any metric space $Y$ equipped with a $\Gamma$-action
by coarse equivalences. 
It is the nonconnective $K$-theory spectrum attached to a symmetric
monoidal category $\mathcal{B}^{\Gamma}_p (Y)$ whose
objects are functors
$\sigma \colon \EGamma \to
\mathcal{B} (\Gamma, \mathcal{B} (Y,R))$
with the additional condition that the morphisms $\sigma(f)$ are bounded by zero but only as homomorphisms between $R$-modules parametrized over $\Gamma$.

The major feature of $\mathcal{B}^{\Gamma}_p (Y)$ is the relaxed control condition on the morphisms in $\mathcal{B} (\Gamma, \mathcal{B} (Y,R))$.  They are still $R$-linear homomorphisms that are bounded in the $\Gamma$-direction.  However the non-zero components between the stalks are allowed to grow in length as the $\Gamma$-coordinates of the stalks are moved further away from the identity in $\Gamma$.
\end{DefRefName}

We note that for metric spaces with isometric $\Gamma$-actions, there is a natural
transformation $K^{\Gamma}_i (Y) \to K^{\Gamma}_p (Y)$.

\begin{ThmRef}{KIOAS}
For the trivial action of $\Gamma$ on $\mathbb{R}^{n}$, there are weak equivalences
\[
\Sigma K^{\Gamma}_i (T\mathbb{R}^{n-1})^{\Gamma} \, \simeq \, K^{\Gamma}_i (T\mathbb{R}^{n})^{\Gamma}
 \, \simeq \, \Sigma^{n+1} \Knc (R\Gamma).
\]
If $\Gamma$ acts on $Y$ by deck transformations with the compact quotient $N$,
then there is a weak equivalence
\[
\sigma \colon \Sigma K^{\Gamma}_i (Y)^{\Gamma} \, \simeq \, K^{\Gamma}_i (T Y)^{\Gamma}.
\]
\end{ThmRef}

\begin{DefRef}{OrbitMet}
Let $Z$ be any metric space with a free left $\Gamma$-action by isometries.
We assume that the action is properly discontinuous, that is, that for fixed points $z$ and $z'$,
the infimum over $\gamma \in \Gamma$ of the distances $d(z, \gamma z')$ is attained.
Then we define the \textit{orbit space metric} on $\Gamma \backslash Z$ by
\[
d_{\Gamma \backslash Z} ([z], [z']) = \inf_{\gamma \in \Gamma} d(z, \gamma z').
\]

Now suppose $X$ is some metric space with a left $\Gamma$-action by isometries.
Define
\[
X^{bdd} = X \times_{\Gamma} \Gamma
\]
where
the right-hand copy of $\Gamma$ denotes $\Gamma$ regarded as a metric space with the word-length metric associated to a finite generating set,
the group $\Gamma$ acts by isometries on the metric space $\Gamma$ via left multiplication,
and $X \times_{\Gamma} \Gamma$ denotes the orbit metric space associated to the diagonal left $\Gamma$-action
on $X \times \Gamma$.
\end{DefRef}

The natural left action of $\Gamma$ on $X^{bdd}$ is given by $\gamma [x,e] = [\gamma x, e]$.
The major feature of this action is that it is always bounded.
Another useful fact is this.

\begin{PropRef}{HJDSEO}
Let $b \colon X \to X^{bdd}$ be the natural map given by $b(x) = [x,e]$ in the orbit space $X \times_{\Gamma} \Gamma$.
The map $b \colon X \to X^{bdd}$ is a coarse map.
\end{PropRef}

If we think of $X^{bdd}$ as the set $X$ with the metric induced from the bijection $b$, the map $b$ becomes the coarse identity map between the metric space $X$ with a left action of $\Gamma$ and the metric space $X^{bdd}$ where the action is made bounded.

One other fact we will need 
is a proper equivariant version of the Spanier-Whitehead duality.

Recall that the compact aspherical manifold $M$ is embedded in a sphere
$S^{n+k}$ for some large $k$, and $Y$ is the universal
cover of the total space $N$ of the normal bundle to
the embedding.
Let $\mathcal{S}$ be a spectrum with an action by $\Gamma$. 
Then there is a weak
equivalence $\Sigma^{n+k} F (X, \mathcal{S})^{\Gamma} \simeq \hlf (Y;
\mathcal{S})^{\Gamma}$. In fact, the canonical map $\mathcal{S}^{\Gamma} \to
\mathcal{S}^{h\Gamma}$ can be identified with a map $\mathcal{S}^{\Gamma} \to
\Omega^{n+k} \hlf (N; \mathcal{S})^{\Gamma}$.

\SSecRef{The beginning of the argument}{Jjjj}

Splitting the map $r_T \colon K(X)^{\Gamma} \to K(X)^{h\Gamma}$ is equivalent to splitting its $(n+k+1)$-fold suspension $\Sigma^{n+k+1} r_T$.
From the fact that $X$ is a model for $E\Gamma$, the codomain can be interpreted as $\Sigma^{n+k+1} F(X, K(X))^{\Gamma}$.
Applying the Spanier-Whitehead duality with $\mathcal{S} = K(X)$, we get a 
weak equivalence
\[
\epsilon \colon 
\Sigma^{n+k+1} F(X, K(X))^{\Gamma} 
\longrightarrow 
\Sigma \hlf (Y; K(X))^{\Gamma}.
\]
This can be followed with
the suspension of the fixed point map
\[
\alpha_{\Gamma} \colon \hlf (Y; K (X))^{\Gamma}
\longrightarrow
K (X \times Y)^{\Gamma}
\]
induced from the equivariant twisted controlled assembly
\[
\phi_{\Gamma} \colon \hlf (Y; K (X))
\longrightarrow
K (X \times Y).
\]
Next we would like to use something like a relax control map $K (X \times Y)^{\Gamma} \to K^{\Gamma}_p (Y)^{\Gamma}$.  However, the free action on $Y$ is not necessarily by coarse equivalences.  
So instead we construct 
\[
\rho \colon \Sigma K (X \times Y)^{\Gamma} \simeq K(X, TY)^{\Gamma} 
\longrightarrow  K^{\Gamma}_p ((TY)^{bdd})^{\Gamma}.
\]
We have designed
a map $t$ as the composition
\[
\rho \circ \Sigma \alpha_{\Gamma} \circ \epsilon \colon
\Sigma^{n+k+1} K(X)^{h\Gamma}
\longrightarrow 
K^{\Gamma}_p ((TY)^{bdd})^{\Gamma}.
\]

If we could compute $K^{\Gamma}_p ((TY)^{bdd})^{\Gamma}$ and show that the composition starting with 
$\Sigma^{n+k+1}  K(X)^{\Gamma}$ is an equivalence, we would be done.
But we can't, and that is likely not true.
So we plan to continue building a map $t$ out to some $\underline{M}$ that is computable and such that the computation results in $\Sigma^{n+k+1}  K(X)^{\Gamma}$
in a particularly natural way.
We do this by introducing a new controlled theory called bounded G-theory developed in \cite{gCbG:00, gCbG:14}. 
This theory has excision properties that are not available in bounded K-theory.
We will review only the necessary details in the next section and then return to the proof in section \refSS{HIUSJ}.

\SSecRef{Fibrewise bounded G-theory}{KICVB}

Let $X$ be a proper metric space and let $R$ be a noetherian ring.
At the basic level bounded G-theory is locally modeled on finitely generated $R$-modules where the exact sequences are not necessarily split.

We will use the notation $\mathcal{P}(X)$ for the power set of $X$ partially ordered by inclusion. Let $\Mod (R)$ denote the category of left $R$-modules.
If $F$ is a left $R$-module, let $\mathcal{I}(F)$ denote the family of all $R$-submodules of $F$ partially ordered by inclusion.
An $X$\textit{-filtered} $R$-module is a module $F$ together with a functor
$\mathcal{P}(X) \to \mathcal{I}(F)$
such that
the value on $X$ is $F$.
It is \textit{reduced} if $F(\emptyset)=0$.

An $R$-homomorphism $f \colon F \to G$ of $X$-filtered modules is \textit{boundedly controlled} if there is a fixed number $b \ge 0$ such that the
image $f (F (S))$ is a submodule of $G (S [b])$
for all subsets $S$ of $X$.
It is called \textit{boundedly
bicontrolled} if, for some fixed $b \ge 0$, in addition to inclusions of submodules
$f (F (S)) \subset G (S [b])$,
there are inclusions
$f (F)
\cap G (S) \subset  f F (S[b])$
for all subsets $S \subset X$.

There is a notion of an admissible exact sequence of filtered modules
which is an exact sequence of $R$-modules with the monic required to be a boundedly bicontrolled monomorphism and the epi required to be a boundedly bicontrolled epimorphism.

\begin{DefRef}{HYUT}
Let $F$ be an $X$-filtered $R$-module.
\begin{enumerate}
\item $F$ is called \textit{lean} or $D$-\textit{lean} if there is a number $D \ge 0$ such that
\[
F(S) \subset \sum_{x \in S} F(x[D])
\]
for every subset $S$ of $X$.
\item $F$ is called \textit{insular} or $d$-\textit{insular} if there is a
number $d \ge 0$ such that
\[
F(S) \cap F(U) \subset F(S[d] \cap U[d])
\]
for every pair of subsets $S$, $U$ of $X$.
\item $F$ is called \textit{split} or $D'$-\textit{split} if there is a
number $D' \ge 0$ such that we have
\[
F(S) \subset F(T[D']) + F(U[D'])
\]
whenever a subset $S$ of $X$ is written as a union $T \cup U$.
\end{enumerate}
\end{DefRef}

Property (3) is a consequence of (1).

It turns out the category $\SI (X)$ of all split insular reduced objects and boundedly controlled morphisms is exact with respect to the family of admissible exact sequences.
The same is true for the subcategory $\LI (X)$ of all lean insular reduced objects.
Both of these categories have cokernels but not kernels.

An $X$-filtered $R$-module $F$ is \textit{locally finitely generated} if $F (S)$ is a finitely generated $R$-module for
every bounded subset $S \subset X$.
The category $\BSI (X)$ is the full
subcategory of $\SI (X)$ on the locally finitely generated objects.
Similarly, the category $\BLI (X)$ is the full
subcategory of $\LI (X)$ on the locally finitely generated objects.
They are closed under extensions, so they are also exact categories.
There is actually a sequence of exact inclusions
$\mathcal{B}(X) \to \BLI (X) \to \BSI (X)$.

There is a delicate issue of existence of so-called gradings of objects in these categories.  We will suppress this issue here and pretend that all $F$ are graded by the submodules $F(S)$.  This is true for objects of $\mathcal{B}(X)$ for example.
In general, one simply requires an appropriate grading.

We now wish to construct a larger \textit{fibred bounded category} $\BGamma (Y)$ which will extend $\mathcal{B}_X (Y) = \mathcal{B} (X, \mathcal{B}(Y))$ similarly to the extension of $\mathcal{B} (X)$ by $\BLI (X)$.
The result will in fact have a mix of features from $\BLI (X)$ and $\BSI (Y)$.

\begin{DefRef}{GEPCW2}
Given an $R$-module $F$, an $(X,Y)$\textit{-filtration} of $F$ is a reduced functor
$\phi_F \colon \mathcal{P}(X \times Y) \to  \mathcal{I}(F)$. 
When there is no ambiguity, we will denote the values $\phi_F (U)$ by $F(U)$.
The associated $X$-filtered $R$-module $F_X$ is given by 
$F_X (S) = F(S \times Y)$.  
Similarly, for each subset $S \subset X$, one has the $Y$-filtered $R$-module $F^S$ given by 
$F^S (T) = F(S \times T)$. 
\end{DefRef}

We will use the following notation generalizing enlargements in a metric space.
Given a subset $U$ of $X \times Y$ and a function $k \colon X \to [0, + \infty )$, let
\[
U [k] = \{ (x,y) \in X \times Y \ \vert \ \textrm{there\ is} \ (x,y') \in U \ \textrm{with} \ d(y,y') \le k(x) \}.
\]
If in addition we are given a number $K \ge 0$ then
\[
U [K,k] = \{ (x,y) \in X \times Y \ \vert \ \textrm{there\ is} \ (x',y) \in U[k] \ \textrm{with} \ d(x,x') \le K  \}.
\]
If $U$ is a product set $S \times T$, it will be convenient to use the notation $(S,T)[K,k]$ in place of $(S \times T)[K,k]$.
We will refer to the pair $(K,k)$ in the notation $U[K,k]$ as the \textit{enlargement data}.

Let $x_0$ be a chosen fixed point in $X$.
Given a monotone function $h \colon [0, + \infty ) \to [0, + \infty )$, there is a function $h_{x_0} \colon X \to [0, + \infty )$ defined by
\[
h_{x_0} (x) = h (d_X (x_0,x)).
\]

Given two $(X,Y)$-filtered modules $F$ and $G$, an $R$-homomorphism $f \colon F(X \times Y) \to G(X \times Y)$ is \textit{boundedly controlled} if
there are a number $b \ge 0$ and a monotone function $\theta \colon [0, + \infty ) \to [0, + \infty )$ such that
\begin{equation}
fF(U) \subset G(U[b,\theta_{x_0}]) \tag{$\dagger$}
\end{equation}
for all subsets $U \subset X \times Y$ and some choice of $x_0 \in X$.
It is clear that the condition is independent of the choice of $x_0$.

\begin{DefRef}{StrFib}
An $(X,Y)$-filtered module $F$ is called
\begin{itemize}
\item \textit{split} or $(D',\Delta')$-\textit{split} if there is a number $D' \ge 0$ and a monotone function
$\Delta' \colon [0, + \infty ) \to [0, + \infty )$ so that
    \[
    F(U_1 \cup U_2) \subset F(U_1 [D',\Delta'_{x_0}]) + F(U_2 [D',\Delta'_{x_0}])
    \]
    for each pair of subsets $U_1$ and $U_2$ of $X \times Y$, \medskip
\item \textit{lean/split} or $(D,\Delta')$-\textit{lean/split} if there is a number $D \ge 0$ and a monotone function
$\Delta' \colon [0, + \infty ) \to [0, + \infty )$ so that \medskip
\begin{itemize}
\item the $X$-filtered module $F_X$ is $D$-lean, while
\item the $(X,Y)$-filtered module $F$ is $(D,\Delta')$-split,
\end{itemize} \medskip
\item \textit{insular} or $(d,\delta)$-\textit{insular} if there is a number $d \ge 0$ and a monotone function
$\delta \colon [0, + \infty ) \to [0, + \infty )$ so that
    \[
    F(U_1) \cap F(U_2) \subset F \big( U_1[d, \delta_{x_0}] \cap U_2[d, \delta_{x_0}] \big)
    \]
    for each pair of subsets $U_1$ and $U_2$ of $X \times Y$.
\end{itemize}
\end{DefRef}

There are two nested categories of $(X,Y)$-filtered modules.
$\mathbf{LS}_X (Y)$ has objects that are lean/split and insular $(X,Y)$-filtered $R$-modules, the morphisms are the boundedly controlled homomorphisms. 
$\B_X (Y)$ is the full subcategory of $\mathbf{LS}_X (Y)$ on objects $F$ such that $F(U)$ is a finitely generated submodule whenever $U \subset X \times Y$ is bounded.  Equivalently, the subcategory $\B_X (Y)$ is full on objects $F$ such that all $Y$-filtered modules $F^S$ associated to bounded subsets $S \subset X$ are locally finitely generated.

A morphism $f \colon F \to G$ in $\U_X (Y)$ is \textit{boundedly bicontrolled} if
there is filtration data $b \le 0$ and $\theta \colon [0, + \infty ) \to [0, + \infty )$, and
in addition to ($\dagger$) one also has the containments
\[
fF \cap G(U) \subset fF(U[b,\theta_{x_0}]).
\]

The category $\B_X (Y)$ is an exact category where the admissible monomorphisms and epimorphisms have to be bicontrolled just as in $\B(X)$.

Suppose $C$ is a subset of $Y$.
Let $\B_X (Y)_{<C}$ be the full subcategory of $\B_X (Y)$ on objects $F$ such that
    \[
    F (X, Y) \subset F \big( (X, C)[r,\rho_{x_0}]) \big)
    \]
    for some number $r \ge 0$ and an order preserving function $\rho \colon [0,+\infty) \to [0,+\infty)$.
It can shown that $\B_X (Y)_{<C}$ is a Grothendieck subcategory of $\B_X (Y)$,
so it is an exact subcategory.
Moreover, there is a homotopy fibration
\begin{gather*}
G_X (Y)_{<C} \longrightarrow G_X (Y) \longrightarrow G_X (Y,C)
\end{gather*}
where $G_X (Y,C)$ is the K-theory of a quotient category $\B_X (Y)/\B_X (Y)_{<C}$.
This localization theorem allows to construct various nonconnective deloopings and is the basic technical tool for proving excision theorems.
The details can be found in \cite{gCbG:13}.
In this paper we take for granted the nonconnective deloopings and simply state the excision theorem.

Given a proper metric space $Y$ with a left $\Gamma$-action by coarse equivalences,
there is the analogue of $K_p^{\Gamma}$ from Definition \refD{POI1}.

\begin{DefRefName}{POI2}{Coarse Equivariant Theories, continued}
Suppose $\Gamma$ acts on $Y$ by coarse equivalences.

$G^{\Gamma}_p (Y)$ is defined as the nonconnective $K$-theory spectrum of the exact category  of functors
$\sigma \colon \EGamma \to \B_{\Gamma} (Y)$
such that the morphisms $\sigma(f)$
are of filtration $0$ when viewed as homomorphisms of $\Gamma$-filtered modules.
We call this category $\B^{\Gamma}_p (Y)$.
\end{DefRefName}

The natural inclusion $\mathcal{B}^{\Gamma}_i (Y) \to \B^{\Gamma}_p (Y)$ gives the \textit{Cartan map} on the fixed points
\[
\kappa^{\Gamma} \colon K_i^{\Gamma} (Y)^{\Gamma} \longrightarrow G_p^{\Gamma} (Y)^{\Gamma}.
\]

Let $Y'$ be a subset of $Y$.
We can relativize the equivariant constructions.

\begin{DefRef}{Nmnnf}
$\B_{\Gamma} (Y)_{<Y'}$ denotes the full subcategory of $\B_{\Gamma} (Y)$ on objects $F$ such that there is a number $k \ge 0$ and an order preserving function $\lambda \colon \mathcal{B}(\Gamma) \to [0,+\infty)$ such that
$F(S) \subset F(S[k])(C[\lambda(S)])$
for each bounded subset $S \subset \Gamma$.  This is a right filtering Grothendieck subcategory.
In particular, there is an exact quotient category $\B_{\Gamma} (Y) / \B_{\Gamma} (Y)_{<Y'}$ which we denote by
$\B_{\Gamma} (Y, Y')$.
\end{DefRef}

Given a proper metric space $Y$ with a left $\Gamma$-action, a subset $Y'$ is called \textit{coarsely $\Gamma$-invariant} or simply \textit{coarsely invariant} if for each element $\gamma$ of $\Gamma$ there is a number $t(\gamma)$ with
$\gamma \cdot Y' \subset Y' [t(\gamma)]$.

\begin{PropRef}{HUIMN}
If $Y'$ is a coarsely invariant subset of $Y$, the subcategory $\B_{\Gamma} (Y)_{<Y'}$ is invariant under the action of $\Gamma$ on $\B_{\Gamma} (Y)$, so there is a left $\Gamma$-action on the quotient $\B_{\Gamma} (Y) / \B_{\Gamma} (Y)_{<Y'}$, and
one obtains the equivariant relative theory $G^{\Gamma}_p (Y,Y')$ and the corresponding fixed point spectrum $G^{\Gamma}_p (Y,Y')^{\Gamma}$.
\end{PropRef}

The construction and the formal properties of the exact category $\W^{\Gamma} (Y,Y')^{\bullet}$ resemble the category of fixed objects $\B_p^{\Gamma} (Y,Y')^{\Gamma}$ while the action is not uniquely specified.

First, we recapitulate the main points in the construction of $\B_p^{\Gamma} (Y,Y')^{\Gamma}$ in a revisionist way that can be generalized.
Suppose $\Gamma$ acts on the proper metric space $Y$ by coarse equivalences and $Y'$ is a coarsely invariant subspace.
The objects of $\B_p^{\Gamma} (Y,Y')$ are the functors $\theta \colon \EGamma \to \B_{\Gamma} (Y,Y')$ such that the morphisms $\theta(f)$ are of filtration $0$ when viewed as homomorphisms of $\Gamma$-filtered modules.
Recall that this category has the left action by $\Gamma$ induced from the diagonal action on $\Gamma \times Y$.
This is the category used in the description of $\B_p^{\Gamma} (Y,Y')^{\Gamma}$.
So, again, a fixed object in $\B_p^{\Gamma} (Y,Y')^{\Gamma}$ is determined by an object
$F$ of $\B_{\Gamma} (Y,Y')$ and isomorphisms $\psi (\gamma) \colon F \to {\gamma} F$ which are of filtration $0$ when projected to $\Gamma$.
We will exploit the fact that this category and its exact structure can be described independently from the construction of the equivariant functor category $\B_p^{\Gamma} (Y,Y')$.
The spectrum $G_p^{\Gamma} (Y,Y')^{\Gamma}$ can be defined as the 
nonconnective K-theory spectrum of $\B_p^{\Gamma} (Y,Y')^{\Gamma}$.

The details of the defininition of $\W^{\Gamma} (Y,Y')^{\bullet}$ and the corresponding nonconnective spectrum $W^{\Gamma} (Y,Y')^{\bullet}$ are mostly straightforward extensions of the summary above.
It should be helpful to point out that the category $\W^{\Gamma} (Y,Y')^{\bullet}$ by itself is not a lax limit with respect to a specific action of $\Gamma$.  However lax limits $\B_p^{\Gamma} (Y,Y')^{\Gamma}$ are going to be exact subcategories of $\W^{\Gamma} (Y,Y')^{\bullet}$.
In order to indicate the lax limit origin of the construction without committing to a specific action, we use the $^{\bullet}$ superscript.

\begin{DefRef}{WGpre2}
The category $\W^{\Gamma} (Y,Y')^{\bullet}$ has objects which are sets of data $( \{ F_{\gamma} \}, \{ \psi_{\gamma} \} )$ where
\begin{itemize}
\item $F_{\gamma}$ is an object of $\B_{\Gamma} (Y,Y')$ for each $\gamma$ in $\Gamma$,
\item $\psi_{\gamma}$ is an isomorphism $F_e \to F_{\gamma}$,
\item $\psi_{\gamma}$ has filtration $0$ when viewed as a morphism over $\Gamma$,
\item $\psi_e = \id$,
\item $\psi_{\gamma_1 \gamma_2} = \gamma_1 \psi_{\gamma_2} \psi_{\gamma_1}$.
\end{itemize}
The morphisms $( \{ F_{\gamma} \}, \{ \psi_{\gamma} \} ) \to ( \{ F'_{\gamma} \}, \{ \psi'_{\gamma} \} )$
are collections $ \{ \phi_{\gamma} \} $, where each $\phi_{\gamma}$ is a morphism $F_{\gamma} \to F'_{\gamma}$ in $\B_{\Gamma} (Y,Y')$, such that the squares
\[
\xymatrix{
F_e \ar[rr]^-{\psi_{\gamma}} \ar[dd]_-{\phi_e} &&F_{\gamma} \ar[dd]^-{\phi_{\gamma}}\\
\\
F'_e \ar[rr]^-{\psi'_{\gamma}} &&F'_{\gamma}
}
\]
commute for all $\gamma \in \Gamma$.
\end{DefRef}

The exact structure on $\W^{\Gamma} (Y,Y')^{\bullet}$ is induced from that on $\B_{\Gamma} (Y,Y')$.
First we observe that for any action $\alpha$ on $Y$ by bounded coarse equivalences, a subset $Y'$ is coarsely invariant, so there is the induced action on the pair $(Y, Y')$.  Now the lax limit category
$\B_p^{\Gamma} (Y,Y')^{\Gamma}_{\alpha}$ is a subcategory of $\W^{\Gamma} (Y,Y')^{\bullet}$.
The embedding $E_{\alpha}$ is realized by sending the object $(F, \psi )$ of $\B_{\Gamma} (Y,Y')^{\Gamma}_{\alpha}$ to
$( \{ \alpha_{\gamma} F \}, \alpha, \{ \psi (\gamma) \} )$.
On the morphisms, $E_{\alpha} (\phi) = \{ {\alpha}_{\gamma} \phi \}$.

A morphism $\phi$ in $\W^{\Gamma} (Y,Y')^{\bullet}$ is an admissible monomorphism if $\phi_e \colon F \to F'$ is an admissible monomorphism in $\B_{\Gamma} (Y,Y')$.  This of course implies that all structure maps $\phi_{\gamma}$ are admissible monomorphisms.  Similarly, a morphism $\phi$ is an admissible epimorphism if $\phi_e \colon F \to F'$ is an admissible epimorphism.
It is clear that for each action $\alpha$ by bounded coarse equivalences the inclusion functor $E_{\alpha}$ is exact.

\begin{DefRef}{JKvvv}
The spectrum $W^{\Gamma} (Y,Y')^{\bullet}$ is the nonconnective K-theory spectrum of $\W^{\Gamma} (Y,Y')^{\bullet}$.
\end{DefRef}

The exact inclusions induce a map of nonconnective spectra
\[
\varepsilon_{\alpha} \colon G_p^{\Gamma} (Y,Y')^{\Gamma}_{\alpha} \longrightarrow W^{\Gamma} (Y,Y')^{\bullet}.
\]
Suppose $\mathcal{U}$ is a finite coarse covering of $Y$.
We define the homotopy colimit
\[
\mathcal{W}^{\Gamma} (Y,Y')^{\bullet}_{< \mathcal{U}} \, = \, \hocolim{U_i \in \mathcal{U}} W^{\Gamma} (Y, Y')^{\bullet}_{<U_i}.
\]
The following excision result is proved in \cite{gCbG:13}.

\begin{ThmRef}{PPPOI8fin}
Suppose the action of $\Gamma$ on $Y$ is by bounded coarse equivalences.
If $\mathcal{U}$ is a finite coarse covering of $Y$ such that the family of all subsets $U$ in $\mathcal{U}$ together with $Y'$ are pairwise coarsely antithetic,
then the natural map induced by inclusions
\[
\delta \colon \mathcal{W}^{\Gamma} (Y,Y')^{\bullet}_{< \mathcal{U}}  \longrightarrow
{W}^{\Gamma} (Y,Y')^{\bullet}
\]
is a weak equivalence.
\end{ThmRef}

\SSecRef{The end of the argument}{HIUSJ}

We paused the construction of $t$ at $K^{\Gamma}_p ((TY)^{bdd})^{\Gamma}$ at the end of section \refSS{Jjjj}.
We now append that map with the following sequence.
The map
\[
\tau \colon
K^{\Gamma}_p ((TY)^{bdd})^{\Gamma}
\longrightarrow
G^{\Gamma}_p ((TY)^{bdd})^{\Gamma}
\]
is induced by the inclusion of equivariant categories.
The map
\[
\mu \colon
G^{\Gamma}_p ((TY)^{bdd})^{\Gamma}
\longrightarrow
W^{\Gamma} ((TY)^{bdd})^{\bullet}
\]
is induced by the inclusion of lax limit categories.

Now we perform a new geometric construction in the normal bundle $N$.
Choose a closed Euclidean ball $B$ in $N$ and
let $\widehat{B}$ be an arbitrary lift of $B$ in $Y$ expressed by a homeomorphism $\sigma \colon B \to \widehat{B}$.
We define the category
\[
\W^{\Gamma}_{\! T \widehat{B}} ((TY)^{bdd}, (T \partial Y)^{bdd})^{\bullet}
=
\W^{\Gamma} ((TY)^{bdd}, (T \partial Y)^{bdd})_{> \complement T \widehat{B}}^{\bullet}
\]
as the exact quotient determined by the complement $\complement T \widehat{B}$ of
$T \widehat{B}$ in $(TY)^{bdd}$.
We have the resulting nonconnective spectrum
\[
W^{\Gamma}_{T \widehat{B}} ((TY)^{bdd}, (T \partial Y)^{bdd})^{\bullet}
=
W^{\Gamma} ((TY)^{bdd}, (T \partial Y)^{bdd})_{> \complement T \widehat{B}}^{\bullet}.
\]
The quotient maps induce the map of spectra
\[
q \colon
W^{\Gamma} ((TY)^{bdd})^{\bullet}
\longrightarrow
W^{\Gamma}_{T \widehat{B}} ((TY)^{bdd}, (T \partial Y)^{bdd})^{\bullet}.
\]

We have finally arrived at the actual target we wished to design for a map $t$ 
in order to split the suspension $\Sigma^{n+k+1} r_T$.
The map $t$ is the composition 
\[
t = q \circ \mu \circ \tau \circ \rho \circ \Sigma \alpha_{\Gamma} \circ \epsilon
\colon
K(X)^{h\Gamma}
\longrightarrow
W^{\Gamma}_{T \widehat{B}} ((TY)^{bdd}, (T \partial Y)^{bdd})^{\bullet}.
\]

\newpage

The way these maps fit together is summarized in the following commutative diagram.

\resizebox{!}{7.5in}{\begin{sideways}
\xymatrix{
%row1
&\Sigma^{n+k+1} K (X)^{\Gamma}
\ar@/^1pc/[r]^{\Sigma^{n+k+1} r_T}
\ar[dddd]^{\simeq}
&\Sigma^{n+k+1} K (X)^{h\Gamma}
\ar@/^1.5pc/[ddr]^{t}
\\
\\
%row2
\Sigma^{n+k+1} K (R\Gamma)
\ar@/^1pc/[uur]^{\simeq}
\ar@/_1pc/[ddr]^{\simeq}
\ar@/_15pc/[rrrr]_{\Sigma^{n+k+1} \kappa}^{\simeq}
&&&W^{\Gamma}_{T \widehat{B}} ((TY)^{bdd}, (T \partial Y)^{bdd})^{\bullet}
\ar^{\ \ \ \ \ e}_{\ \ \ \ \ \simeq}[r]
&\Sigma^{n+k+1} G (R\Gamma) \ \ 
\\
\\
%row4
&K^{\Gamma}_i (T\mathbb{R}^{n+k})^{\Gamma}
\ar[r]
&K^{\Gamma}_i (T{N},T \partial {N})^{\Gamma}
\ar@/_1.5pc/[uur]^{\lambda}
}
\end{sideways}
}

\newpage

It remains to apply the excision theorem to check that there is a weak equivalence
\[
e \colon
W^{\Gamma}_{T \widehat{B}} ((TY)^{bdd}, (T \partial Y)^{bdd})^{\bullet}
\longrightarrow
\Sigma^{n+k+1} G (X)^{\Gamma}.
\]
Once that is done, the two ingredients in needed to complete the proof of injectivity of $r_T$
are the facts 
\begin{enumerate}
\item that 
\[
e \circ \phi \circ \Sigma^{n+k+1} \rho
\colon
\Sigma^{n+k+1} K (X)^{\Gamma}
\longrightarrow
\Sigma^{n+k+1} G (X)^{\Gamma}
\]
is homotopic to the $(n+k+1)$-fold suspension of the Cartan map and 
\item that the Cartan map \[ \kappa \colon K (R\Gamma) \longrightarrow G (R\Gamma) \] is an equivalence under the assumptions of Theorem \refT{MTBB}.
\end{enumerate}
Verifying the fact (1) requires a major effort in \cite{gCbG:13}.
It is based on constructing and analizing the infinite transfer map 
\[
\lambda \colon 
K^{\Gamma}_i (T{N},T \partial {N})^{\Gamma} 
\longrightarrow 
W^{\Gamma}_{T \widehat{B}} ((TY)^{bdd}, (T \partial Y)^{bdd})^{\bullet}.
\]
The fact (2) is the weak regular coherence property of $\Gamma$ that was proved by the authors in \cite{gCbG:03,gCbG:13}.
Remarkably, the class of groups which are weakly regular coherent turns out to contain all of the groups for which the Novikov conjecture has been verified.  Since the work of Waldhausen \cite{fW:78} variants of this property were thought to be mysterious and somewhat a bottleneck in bridging proofs of injectivity to bijectivity for assembly maps.  This turns out not the case within this approach.  In the last section \refSS{CCFG}, we will describe how far this approach is expected to extend the class of groups and rings for which the analogue of Theorem \refT{MTBB} can be verified.

It remains to construct the map $e$.
That will require specific subsets of $T Y$.

We may assume that the chosen ball $B$ is a metric ball in $\mathbb{R}^{n+k}$ centered at $c$ with radius $R$ and contained entirely in the interior of $N$.
Let $h \colon \mathbb{R}^{n+k} \to \mathbb{R}^{n+k}$ be the linear map $h(x) = c + Rx$,
so $h$ restricts to a linear homeomorphism
$h \colon D^{n+k} \to B$ from the unit disk $D^{n+k} = 0[1]$ onto the chosen ball $B \subset N$.

\begin{NotRef}{OPIEWQp}
We identify the following subsets of $\mathbb{R}^{n+k}$:
\begin{align}
&E^{+}_i = \{ (x_1, \ldots, x_{n+k}) \mid x_l = 0 \mathrm{ \ for\ all \ } l > i,\ x_i \ge 0 \}, \notag \\
&E^{-}_i = \{ (x_1, \ldots, x_{n+k}) \mid x_l = 0 \mathrm{ \ for\ all \ } l > i,\ x_i \le 0 \}, \notag \\
&E_i = E^{-}_i \cup E^{+}_i, \mathrm{ \ for\ } 1 \le i \le n+k, \mathrm{ \ and\ } \notag \\
&E_0 = E^{-}_0 = E^{+}_0 = \{ (0, \ldots, 0) \}. \notag
\end{align}
There are also related to $E^{\pm}_i$ subsets
\begin{align}
&D_0 = D^{-}_0 = D^{+}_0 = \{ (0, \ldots, 0) \}, \notag \\
&D^{\pm}_i = E^{\pm}_i \cap D^{n+k}, \mathrm{ \ for\ } 1 \le i \le n+k, \mathrm{ \ and\ } \notag \\
&D_i = D^{-}_i \cup D^{+}_i = D^{-}_{i+1} \cap D^{+}_{i+1}. \notag
\end{align}
The images of $D^{\ast}_i$ under the linear homeomorphism $h$ will be called $B^{\ast}_i$.
\end{NotRef}

The subsets $E^{\ast}_i$ form a coarsely antithetic covering of $\mathbb{R}^{n+k}$.
It is easy to see that $TE^{\ast}_i$ form a coarsely antithetic covering of $T\mathbb{R}^{n+k}$.
Therefore, we obtain a coarsely antithetic covering $\mathcal{E}$ of $T\mathbb{R}^{n+k}$ by the subsets $T h(E^{\ast}_i)$ .
Notice that this covering is closed under coarse intersections since it includes the subsets $T h(E_i)$ for
$0 \le i < n+k$, where $E_i$ are the intersections $E^{-}_{i+1} \cap E^{+}_{i+1}$.

Now we define the following collection of subsets of $V = T \widehat{B}$:
\[
V' = T \partial \widehat{B}, \ V^{\ast}_i = T \sigma (B^{\ast}_i ), \ \textrm{for}\ 0 \le i \le n+k,\ \textrm{and}\ V_i = \left( T \sigma h (D_i) \right). 
\]
The subsets $\{ V^{\ast}_i \}$ can be thought of as a coarsely antithetic covering of $T \widehat{B}$.
They can be extended to a coarsely antithetic covering of the metric space $\overline{V} = (T Y )^{bdd}$:
\[
\overline{V'} = T \partial \widehat{B} \cup \complement T \widehat{B}, \ 
\overline{V^{\ast}_i} = T \sigma (B^{\ast}_i ) \cup \complement T \widehat{B}  \mathrm{ \ \ for \ } 0 \le i \le n+k.
\]

Let $\mathcal{U}$ be the coarse antithetic covering of $T \widehat{B}$ by $V^{\ast}_i$.
There is a homotopy pushout
\[
\mathcal{W}^{\Gamma}_{T \widehat{B}} ((TY)^{bdd}, (T \partial Y)^{bdd})_{<  \mathcal{U}} =
\hocolim{\mathcal{U}} W^{\Gamma}_{T \widehat{B}} ((TY)^{bdd}, (T \partial Y)^{bdd})_{< \overline{V^{\ast}_i}}^{\bullet}.
\]

From Theorem \refT{PPPOI8fin} we have the following consequence.

\begin{ThmRef}{njuyhb}
There is a weak equivalence
\[
 \delta \colon \mathcal{W}^{\Gamma}_{T \widehat{B}} ((TY)^{bdd}, (T \partial Y)^{bdd})_{< \mathcal{U}}^{\bullet}
\longrightarrow
{W}^{\Gamma} ((TY)^{bdd}, (T \partial Y)^{bdd})^{\bullet}.
\]
So we get a weak equivalence
\[
e \colon
W^{\Gamma}_{T \widehat{B}} ((TY)^{bdd}, (T \partial Y)^{bdd})^{\bullet}
\longrightarrow
\Sigma^{n+k+1} G (R\Gamma).
\]
\end{ThmRef}

In \cite{gCbG:13,gCbG:15} we establish that the composition
\[
\Sigma^{n+k+1} K (R\Gamma) \to \Sigma^{n+k+1} G (R\Gamma)
\]
is a map of the homotopy colimits where all component maps are equivalences.
Therefore, the composite is indeed the suspension  $\Sigma^{n+k+1} \kappa$ and is an equivalence according to fact (2).

\SSecRef{Coarse coherence of fundamental groups}{CCFG}

In this section we survey the resolution of fact (2) above under the name of \textit{weak regular coherence} in \cite{gCbG:03,gCbG:13,bG:13,bGjG:18}, generalizing the work of Waldhausen \cite{fW:78}. 

For a ring $A$, a \textit{presentation} of an $A$-module $E$ is an exact sequence  $F_2 \to F_1 \to E \to 0$ with both $F_1$ and $F_2$ free $A$-modules. It is a \textit{finite presentation} if the free modules are finitely generated.  More generally, one has the notion of a \textit{projective resolution} of $E$ which is an exact sequence 
\[
\ldots \longrightarrow P_n \longrightarrow \ldots \longrightarrow P_2 \longrightarrow P_1 \longrightarrow E \longrightarrow 0
\]
where all $P_i$ are projective $A$-modules.  The projective resolution is of \textit{finite type} if the projective modules are finitely generated.  It is called \textit{finite} if there is a number $n$ such  that the modules $P_i =0$ for $i > n$.
A module is said to have \textit{finite projective dimension} if
it has a projective resolution of finite type.
In turn, the ring $A$ has \textit{finite global dimension}
if there is a number $n$ such that every finitely generated $A$-module has a finite projective resolution of length $n$. 

The ring $A$ is called \textit{coherent} if every finitely presented $A$-module has finite projective dimension.
A coherent ring $A$ is called \textit{regular coherent} if each projective resolution of finite type over $A$ is chain homotopy equivalent
to a finite projective resolution.
Restricting further, a \textit{regular Noetherian} ring is a regular coherent ring which is Noetherian.  Specializing to groups, Waldhausen called a group $\Gamma$ \textit{regular coherent} if the group algebra $R\Gamma$ is
regular coherent for any choice of a regular Noetherian ring $R$.

The collection $\frak{X}$ of regular coherent groups includes free groups, free abelian groups,
torsion-free one relator groups, fundamental groups of
submanifolds of the three-dimensional sphere, and their various
amalgamated products and HNN extensions and so, in particular, the
fundamental groups of submanifolds of the three-dimensional sphere.  Waldhausen used this property to compute the algebraic $K$-theory of
regular coherent groups. 

Two remarks regarding Waldhausen's regular coherence are in order.

(1) The regular coherence property seems to be very special: simply constructing
individual non-projective finite dimensional modules over group rings is hard.  

(2) The collection $\frak{X}$ is not well-understood structurally beyond the portion identified by Waldhausen.  For example, it is unknown whether $\frak{X}$ is closed under products.  While all groups in $\frak{X}$ are necessarily torsion-free, it is unknown if there is a torsion-free group outside of it.

Waldhausen asked if a weaker property of the group
ring would suffice in his argument,
see for example the paragraph after the proof of Theorem 11.2 in~\cite{fW:78}.
A resolution of this question from \cite{bGjG:18} needs an extra condition in the spirit of Definition \refD{HYUT}.

An $X$-filtered $R$-module
$F$ is called \textit{scattered} or $\delta$-\textit{scattered} if there is a number $\delta \ge 0$ such that
\[
F(X) \subset \sum_{x \in X} F(x[\delta]).
\]
This is a consequence of the lean property but is the property responsible for finite generation of the modules over $R\Gamma$.  
A metric space $X$ is called \textit{coarsely coherent} if in any exact sequence 
\[
0 \to E' \xrightarrow{\ f \ } E \xrightarrow{\ g \ } E'' \to 0
\]
of $X$-filtered $R$-modules where $f$ and $g$ are both bicontrolled maps, the combination of $E$ being lean and $E''$ being insular implies that $E'$ is necessarily scattered.

It is shown in \cite{bG:13} that a metric space $X$ with the geometric property called straight finite decomposition complexity (sFDC) in the sense of \cite{aDmZ:12} is coarsely coherent. 
 In particular, a metric space $X$ with finite asymptotic dimension is coarsely coherent which was essentially known already from \cite{gCbG:03}.

A finitely generated group $\Gamma$ is \textit{coarsely regular coherent relative to a ring $R$} if every $R[\Gamma]$-module $F$ which is lean, insular, and locally finitely generated when viewed as a $\Gamma$-filtered $R$-module has a finite projective resolution over $R[\Gamma]$.  It is called simply \textit{coarsely regular coherent} if it is coarsely regular coherent relative to any regular coherent ring of finite global dimension.

It turns out that the proof of the main theorem of \cite{gCbG:03} now gives this consequence.

\begin{ThmRef}{SFDCCRC}
	A finitely generated group $\Gamma$ with sFDC and a finite model for $K(\Gamma,1)$ is coarsely regular coherent.
\end{ThmRef}

We will use the notation $\frak{Y}$ for the class of all coarsely regular coherent groups.
The main interest in this class is this fact.

\begin{ThmRef}{SFDCCRC2}
	A finitely generated group $\Gamma$ which is coarsely regular coherent is weakly regular coherent.  Therefore, for any regular Noetherian ring $R$ of finite global dimension, the Cartan map $ \kappa \colon K (R\Gamma) \to G (R\Gamma) $ is an equivalence
\end{ThmRef}

In order to introduce structure into the class $\frak{Y}$, one can use the usual in coarse geometry so-called permanence properties.  Indeed, coarse coherence satisfies all of the most important permanence properties: Fibering Permanence, Subspace Permanence, Finite Amalgamation Permanence, Finite Union Permanence, Union Permanence and Limit Permanence.  We refer to \cite{bGjG:18} for details.  These geometric facts lead to the invariance of $\frak{Y}$ under, among other constructions, extensions such as finite semi-direct products of groups, amalgamated free products, and HNN extensions.
An elementary consequence of this fact is that the family $\frak{X}$ of regular coherent groups of Waldhausen, which are essentially multiple amalgamated free products and HNN extensions of finitely many copies of $\mathbb{Z}$, are in $\frak{Y}$.

\SecRef{Comparison with the Farrell-Jones conjecture}{fjcomp}

The Farrell-Jones isomorphism conjectures were introduced by F.T.~Farrell and L.E.~Jones  in \cite{fFlJ:93} as an attempt to model the K-theory of the group ring in terms of group homology up to the contribution from the virtually cyclic subgroups.
They verified the conjectures for a variety of fundamental groups of manifolds, most notably the nonpositively curved Riemannian manifolds \cite{fFlJ:86}, linear groups \cite{fFlJ:93}, and other classes of fundamental groups \cite{fFpL:86,fFlJ:87,fFlJ:88}.
A more flexible formulation of the conjectures by J. Davis and W. L\"{u}ck in \cite{jDwL:98} is in terms of orbit categories for families of subgroups. There are many excellent expositions of the Farrell-Jones conjecture such as \cite{aBwLhR:08,wLhR:05}.  We specifically recommend the up-to-date survey by Reich and Varisco \cite{hRmV:18}.

A \textit{family of subgroups} $\mathcal{F}$ of a group $\Gamma$ is a nonempty family closed under conjugation and under passage to subgroups. 
Examples of such families are the collections of all finite subgroups $\mathit{Fin}$,
all virtually cyclic subgroups $\mathit{VCyc}$,
and the two extreme cases of the only trivial subgroup $1$ and the family of all subgroups denoted simply by $\Gamma$.

The orbit category of the group $\Gamma$ is the category $\Or_\Gamma$ of all coset spaces of $\Gamma$ regarded 
as $\Gamma$-spaces, so the morphisms between $\Gamma/H$ and $\Gamma/K$ are the $\Gamma$-maps $\Gamma/H \to \Gamma/K$.

Given a family $\mathcal{F}$ of subgroups of $\Gamma$, one defines the full 
subcategory $\Or_\Gamma \mathcal{F}$ on the coset spaces $\Gamma/F$ where $F$ belongs to the family 
$\mathcal{F}$.
So in particular $\Or_\Gamma \Gamma = \Or_\Gamma$.

The canonical maps $K(R) \to K(R\Gamma)$ induce the assembly map 
\[
A(1) \colon \hocolim{\Or_\Gamma 1} K(R) \longrightarrow K(R\Gamma).
\]
One can check that this is precisely the Loday assembly from section \refSS{loday}.

For a general family $\mathcal{F}$, there is a way to define a functor 
$K_R \colon \Or_\Gamma \to \textit{spectra}$ so that
for all subgroups $F < \Gamma$ there is a homotopy equivalence 
$K_R (\Gamma/F) \simeq K(RF)$.
Using this functor and the canonical maps $K(RF) \to K(R\Gamma)$,
one obtains the Farrell-Jones assembly map 
\[
A(\mathcal{F}) \colon \hocolim{\Or_\Gamma \mathcal{F}} K_R \longrightarrow K(R\Gamma).
\]
The Farrell-Jones conjecture relative to the family $\mathcal{F}$ is the statement that the assembly map $A(\mathcal{F})$ is an equivalence.

The following diagram illustrates the relationship between various maps that have appeared in this paper.
\[
\xymatrix{
 \hocolim{\Or_\Gamma 1} G(R)  \ar[rrr]^-{A_G (1)}
&&&G(R\Gamma) \\ 
\\
 \hocolim{\Or_\Gamma 1} K(R) \ar[rrr]^-{A(1)} \ar[dd]_-{\gamma}
 \ar[uu]^-{\hocolimprep \kappa}
&&& K(R\Gamma) \ar[uu]_-{\kappa_\Gamma}
\\
\\
\hocolim{\Or_\Gamma \mathit{VCyc}} K_R 
\ar@/_1.3pc/[rrruu]_-<<<<<{\ A(\mathit{VCyc})}  
}
\]
The assembly map relevant to the study of topological rigidity of closed aspherical manifolds is $A(1)$.  Even more specifically, the assembly map of genuine geometric interest is 
\[
A(1) \colon \hocolim{\Or_\Gamma 1} K(\mathbb{Z}) \longrightarrow K(\mathbb{Z}\Gamma),
\]
where $\Gamma$ is the fundamental group of a closed aspherical manifold, and is therefore torsion-free.

Now we observe the following facts.  The Cartan map $\kappa \colon K(R) \to G(R)$ is an equivalence for the regular noetherian ring $\mathbb{Z}$.  The map $\gamma$ is also an equivalence whenever $R$ is regular noetherian and $\Gamma$ is torsion-free, which is a much more subtle fact.  Finally, for a very large class of groups $\Gamma$ with strict finite decomposition complexity, the Cartan map $\kappa_\Gamma$ on the right is an equivalence.  So, in this geometrically important and very general situation, the two extreme maps $A_G (1)$ and $A (\mathit{VCyc})$ coincide up to homotopy.
In both lines of research the fact is that studying these maps is technically and conceptually easier than $A (1)$ itself.

The two parts of the diagram above and below $A(1)$ are two different approaches to analyzing the map.  The Farrell-Jones conjecture is an attempt to modify the domain of $A(1)$ with the inherited problems of analyzing and possibly adjusting the domain.  The G-theoretic approach is essentially an attempt to modify the target of $A(1)$ with the inherited need to compute the Cartan map.

It should be interesting to compare these two approaches with the goal of identifying the strengths of each method.  

The strength of the work on the Farrell-Jones conjecture, even the latest general advances 
\cite{aBwLhR:08a, aBwL:12, aBwLhRhR:14,aBmB:18}, is in its geometric nature.  The isomorphism results about the Farrell-Jones assembly come with various geometric conditions on the $\Gamma$-spaces but are true for any choice of the ring $R$.  So in order to relate these results to $A (1)$, one simply needs to assume that the coefficient ring $R$ is regular noetherian, which is a very weak assumption.
On the other hand, the contributions in this direction still use (up to some induction procedures) the flows on $\Gamma$-spaces, and therefore some possibly immanent version of nonpositive curvature.
This geometric nature restricts the range of accessible group geometries.

Turning to our approach, the transition from $A (1)$ to the more algebraic map $A_G (1)$, which is an equivalence for any $\Gamma$ with finite $B\Gamma$ and any noetherian $R$, also requires constraints on both $R$ and $\Gamma$.  However, in this method the algebraic constraints on $R$ are more strict (finite homological dimension) but the geometric constraints are considerably more lax and abstract (strict finite decomposition complexity).  The assembly map is addressed directly in this method, not through the study of the cokernel.  This allows us, in all cases when we prove the isomorphism conjecture, exhibit the inverse isomorphism.  Similar to statements of the Farrell-Jones conjecture, it should be possible to incorporate torsion in this method when the coefficient ring is, for example, the rationals $\mathbb{Q}$. 

It would be fantastic to find some middle ground where the two approaches through the G-theory assembly and through the Farrell-Jones assembly can be made to interact and hopefully enrich each other.

One intriguing step in this direction is in the recent work of Corti\~{n}as and Cirone \cite{gCeC:14}.
A corollary to their main theorem states that if the Farrell-Jones conjecture is true for a given group and coefficients in any commutative smooth $\mathbb{Q}$-algebra $R$ then it is also true for the same group and any commutative $\mathbb{Q}$-algebra.
Since every smooth $\mathbb{Q}$-algebra $R$ is noetherian regular of finite homological dimension, the Cartan map $\kappa_\Gamma$ in the diagram is an equivalence for a geometric group $\Gamma$ of finite decomposition complexity.
The fact that the assembly map $A_G (1)$ is an equivalence gives the Farrell-Jones in the special smooth cases, therefore $A(\mathit{VCyc})$ is also an equivalence for these groups and coefficients in any commutative $\mathbb{Q}$-algebra.


\begin{thebibliography}{99}

\bibitem{fA:78}
J.F. Adams,
\textit{Infinite loop spaces}, 
Ann. Math. Studies 90, Princeton University Press, 1978.

\bibitem{fA:95}
\bysame, 
\textit{Stable homotopy and generalised homology}, 
reprint of the 1974 original, Chicago Lectures in Mathematics, University of Chicago Press, 1995.

\bibitem{aB:03}
{A.~Bartels}, \textit{Sqeezing and higher algebraic \textit{K}-theory}, 
K-theory {\bf 28},
(2003), 19--37.

\bibitem{aBmB:18}
A. Bartels and M. Bestvina,
\textit{The Farrell-Jones Conjecture for mapping class groups},
\texttt{arXiv:1606.02844}

\bibitem{aBwLhR:08}
A.~Bartels, W.~L\"{u}ck, and H.~Reich, 
\textit{On the Farrell-Jones Conjecture and its applications},
J. Topol. \textbf{1} (2008), 57--86.

\bibitem{aBwLhR:08a}
\bysame,
\textit{The K-theoretic Farrell-Jones Conjecture for hyperbolic groups},
Invent. Math. (2008), 29--70.

\bibitem{aBwL:12}
A.~Bartels and W.~L\"{u}ck
\textit{The Borel Conjecture for hyperbolic and CAT(0)-groups},
Ann. Math. 175 (2012), 631--689.

\bibitem{aBwLhRhR:14}
A.~Bartels, W.~L\"{u}ck, H.~Reich, and H.~R\"{u}ping,
\textit{K- and L-theory of group rings over $GL_n (\mathbb{Z})$}, 
Publ. Math. Inst. Hautes \'{E}tudes Sci.  \textbf{119} (2014), 97--125.

\bibitem{mBkBkF:10}
{M.~Bestvina, K.~Bromberg, and K.~Fujiwara},
\textit{Constructing group actions on quasi-trees and applications to mapping class groups},
Publ. Math. Inst. Hautes \'{E}tudes Sci. \textbf{122} (2015), 1--64.

\bibitem{mBwHiM:93}
M. B\"{o}kstedt, W.C. Hsiang, and I. Madsen, 
\textit{The cyclotomic trace and algebraic $K$-theory of spaces}, Invent. Math. \textbf{111} (1993), 465--539.

\bibitem{aBdK:72}
{A. K. Bousfield and D. M. Kan},
{\it Homotopy limits, completions, and localizations},
Lecture Notes in Mathematics {\bf 304}, Springer-Verlag (1972).

\bibitem{gC:92}
{G. Carlsson},
{\it Equivariant stable homotopy theory},
Topology \textbf{31} (1992), 1--27.

\bibitem{gC:93}
\bysame, 
{\it Proper homotopy theory and transfers for infinite groups}, 
in Algebraic Topology and its Applications,
MSRI Publications \textbf{27} Springer-Verlag (1993), 1--10.

\bibitem{gC:95}
\bysame, 
{\it Bounded K-theory and the assembly map
in algebraic K-theory},
in {\it Novikov conjectures, index theory and rigidity},
{\it Vol. 2}
(S.C. Ferry, A. Ranicki, and J. Rosenberg, eds.),
Cambridge U. Press (1995), 5--127.

\bibitem{gC:05}
\bysame,
{\it Deloopings in Algebraic K-Theory}, 
in \textit{Handbook of K-theory} (E. Friedlander, D. Grayson, eds.), Springer (2005), 3--38.

\bibitem{gCbG:03}
{G.~Carlsson and B.~Goldfarb}, \textit{On homological coherence of
discrete groups}, J.~Algebra {\bf 276} (2004), 502--514.

\bibitem{gCbG:04}
\bysame, \textit{The integral K-theoretic Novikov conjecture for
groups with finite asymptotic dimension}, Invent. Math. {\bf
157} (2004), 405--418.

\bibitem{gCbG:00}
\bysame,
\textit{Controlled algebraic G-theory, I}, J. Homotopy Relat. Struct. \textbf{6} (2011), 119--159.

\bibitem{gCbG:13}
\bysame, \textit{Algebraic K-theory of geometric groups}, 2016.
\texttt{arXiv:1305.3349}

\bibitem{gCbG:14}
\bysame, \textit{K-theory with fibred control}, 2016.
\texttt{arXiv:1404.5606}

\bibitem{gCbG:15}
\bysame, \textit{On modules over infinite group rings}, Int. J. Algebra Comput. \textbf{26} (2016), 1--16.

\bibitem{gCjM:95}
G. Carlsson and R.J. Milgram, 
\textit{Stable homotopy and iterated loop spaces}, in
\textit{Handbook of Algebraic Topology} (I.M. James, ed.), North-Holland (1995), 505--583.

\bibitem{gCeP:93}
{G. Carlsson and E.K. Pedersen}, {\it Controlled algebra and the
Novikov conjecture for K- and L-theory}, Topology {\bf 34}
(1993), 731--758.

\bibitem{gCeP:98}
\bysame, {\it \v{C}ech homology and the Novikov conjectures for
K- and L-theory}, Math. Scand. \textbf{82} (1998), 5--47.

\bibitem{sC:04}
S. Chang,
\textit{On conjectures of Mathai and Borel}, Geom. Dedicata \textbf{106} (2004), 161--167. 

\bibitem{sCsFgY:08}
S. Chang, S. Ferry, and G. Yu,
\textit{Bounded rigidity of manifolds and asymptotic dimension growth},
J. K-Theory \textbf{1} (2008), 129--144.

\bibitem{tC:74}
T. Chapman, 
\textit{Topological invariance of Whitehead torsion}, 
Ann. Math. \textbf{96} (1974), 488--497.

\bibitem{gCeC:14}
G.~Corti\~{n}as and E.R.~Cirone,
\textit{Singular coefficients in the K-theoretic Farrell-Jones conjecture}, Algebr. Geom. Topol. \textbf{16} (2016), 129--147.
\texttt{arXiv:1403.1855}

\bibitem{jD:00}
J.F. Davis,
\textit{Manifold aspects of the Novikov conjecture},
in \textit{Surveys on surgery theory, Vol. 1},  Princeton U. Press (2000), 195--224.

\bibitem{jDwL:98}
J.F. Davis and W. L\"{u}ck, 
\textit{Spaces over a category and assembly maps in isomorphism conjectures in K-and L-theory}, 
K-Theory \textbf{15} (1998), 201--252.

\bibitem{aD:05}
A. Dranishnikov,
\textit{Dimension Theory: Local and Global},
in Proceedings of the 19th Annual Workshop in Geometric Topology,
Calvin College, Grand Rapids, MI (June 13-15, 2002).

\bibitem{aD:04}
{\bysame}, 
\textit{Lipschitz cohomology, Novikov's conjecture, and expanders}, Tr. Mat. Inst. Steklova \textbf{247} (2004), 59--73 (Russian, with Russian summary); English transl., Proc. Steklov Inst. Math. \textbf{247} (2004), 50--63.

\bibitem{aDsFsW:08}
A. Dranishnikov, S. Ferry, S. Weinberger,
\textit{An Etale approach to the Novikov conjecture},
Pure Appl. Math. \textbf{61} (2008), 139--155.

\bibitem{aDmZ:12}
{A.~Dranishnikov and M.~Zarichnyi}, \textit{Asymptotic dimension, decomposition complexity, and Haver's property C}, Topology Appl. \textbf{169} (2014), 99--107.

\bibitem{fF:02}
F.T. Farrell,
{\it The Borel conjecture},
in {\it Topology of high-dimensional manifolds (Trieste, 2001)}, ICTP Lecture Notes \textbf{9} (2002), 225--298.

\bibitem{fFlJ:86}
F.T. Farrell and L.E. Jones,
\textit{K-theory and dynamics. I}, 
Ann. Math. \textbf{124} (1986), 531--569.

\bibitem{fFlJ:87}
{\bysame},
\textit{Algebraic K-theory of discrete subgroups of Lie groups}, 
Proc. Nat. Acad. Sci. U.S.A. \textbf{84} (1987), 3095--3096.

\bibitem{fFlJ:88}
{\bysame},
\textit{The surgery L-groups of poly-(finite or cyclic) groups}, Invent. Math. \textbf{91} (1988), 559--586.

\bibitem{fFlJ:91}
{\bysame}, 
\textit{Rigidity in geometry and topology}, 
in Proceedings
of the International Congress of Mathematicians, Vol. I, II (Kyoto, 1990), Math. Soc. Japan. (1991), 653--663.

\bibitem{fFlJ:93}
{\bysame}, 
{\it Isomorphism conjectures in algebraic K-theory}, 
J. Amer. Math. Soc. \textbf{6} (1993), 249--297.

\bibitem{fFlJ:93}
{\bysame}, 
{\it Rigidity for aspherical manifolds with $\pi_1 \subset \mathrm{GL}_m (\mathbb{R})$}, Asian J. Math. \textbf{2} (1998), 215--262.

\bibitem{fFpL:86}
F.T. Farrell and P.A. Linnell,
\textit{K-theory of solvable groups}, 
Proc. London Math. Soc. \textbf{87} (2003), 309--336.

\bibitem{sF:81}
S. Ferry,
\textit{A simple-homotopy approach to the finiteness obstruction},
in \textit{Shape theory and geometric topology (Dubrovnik, 1981)}, Lecture Notes in Mathematics, Springer-Verlag (1981), 73--81.

\bibitem{sFaR:01}
S. Ferry and A. Ranicki,
\textit{A survey of Wall's finiteness obstruction},
in \textit{Surveys on Surgery Theory, Vol. 2}, Ann. Math. Studies 149, Princeton U. Press (2001), 63--80.

\bibitem{OB:95}
S.C. Ferry, A. Ranicki, and J. Rosenberg, eds. 
{\it Novikov conjectures, index theory and rigidity},
{\it Vols. 1, 2},
Cambridge U. Press (1995).

\bibitem{OB:95s}
{\bysame},
\textit{A history and survey of the Novikov Conjecture}, 
in {\it Novikov conjectures, index theory and rigidity},
Vol. 1,
Cambridge U. Press (1995), 7--66.

\bibitem{bG:97}
{B. Goldfarb}, {\it Novikov conjectures for arithmetic groups
with large actions at infinity}, K-theory {\bf 11} (1997),
319--372.

\bibitem{bG:98}
{\bysame}, 
{\it Novikov conjectures and relative hyperbolicity},
Math.~Scand. \textbf{85} (1999), 169--183.

\bibitem{bG:99}
{\bysame}, 
\textit{Large scale topology and algebraic K-theory of arithmetic groups},
Topol. Appl. \textbf{140} (2004), 267--294.

\bibitem{bG:13}
{\bysame}, 
\textit{Weak coherence and the K-theory of groups with finite decomposition complexity}, to appear in Int. Mat. Res. Not. IMRN.
\texttt{arXiv:1307.5345}

\bibitem{bGjG:18}
{B. Goldfarb and J.L. Grossman},
\textit{Coarse coherence of metric spaces and groups and its permanence properties}, Bull. Austr. Math. Soc. \textbf{98} (2018), 422--433.

\bibitem{mG:96}
M. Gromov,
\textit{Positive curvature, macroscopic dimension, spectral gaps and higher signatures}, 
in \textit{Functional analysis on the eve of the 21st century}, Vol. II (New Brunswick, NJ, 1993), Progress in Mathematics 132, 
Birkh\"{a}user (1996), 1--213.

\bibitem{eGnHsW:05}
E. Guentner, N. Higson and S. Weinberger,
\textit{The Novikov conjecture for linear groups},
Publ. Math. Inst. Hautes \'{E}tudes Sci. \textbf{101} (2005), 243--268.

\bibitem{iH:02}
I. Hambleton,
\textit{Algebraic K- and L-theory and applications to topology of manifolds},
in {\it Topology of high-dimensional manifolds (Trieste, 2001)}, ICTP Lecture Notes \textbf{9} (2002), 299--369.

\bibitem{iHeP:04}
I. Hambleton and E.K. Pedersen,
\textit{Identifying assembly maps in K- and L-theory},
Math. Ann. \textbf{328} (2004), 27--57. 

\bibitem{mKwL:05}
M. Kreck and W. L\"{u}ck,
\textit{The Novikov Conjecture: Geometry and Algebra}, 
Birkhauser (2005).

\bibitem{kL:92}
{K. Liu},
{\it On mod 2 and higher elliptic genera},
Comm. Math. Phys. {\bf 149} (1992), 71--95.

\bibitem{jlL:76}
{J.-L.~Loday},
\textit{K-th\'eorie alg\'ebrique et repr\'esentations de groupes},
Ann. sci. \'Ecole Norm. Sup. \textbf{9} (1976), 309--377.

\bibitem{wLhR:05}
W. L\"{u}ck and H. Reich,
\textit{The Baum-Connes and the Farrell-Jones Conjectures in K- and L-Theory},
in \textit{Handbook of K-theory} (E. Friedlander, D. Grayson, eds.), Springer (2005), 703--842.

\bibitem{jpM:77}
J.P. May,
\textit{Infinite loop space theory},
Bull. Amer. Math. Soc. \textbf{83} (1977), 456--494.

\bibitem{gM:95}
G. Mislin,
\textit{Wall's finiteness obstruction}, 
in \textit{Handbook of algebraic topology},
North-Holland (1995), 1259--1291.

\bibitem{sN:64}
S.P. Novikov, 
\textit{Homotopically equivalent smooth manifolds. I} 
Izv. Akad. Nauk SSSR Ser. Mat. \textbf{28} (1964), 365--474.

\bibitem{sN:65a}
{\bysame}, \textit{Rational Pontrjagin classes. Homeomorphism and homotopy type of closed manifolds. I},
Izv. Akad. Nauk SSSR Ser. Mat. \textbf{29} (1965), 1373--1388.

\bibitem{sN:65b}
{\bysame}, \textit{Topological invariance of rational classes of Pontrjagin}, Dokl. Akad. Nauk SSSR \textbf{163} (1965), 298--300.

\bibitem{sN:10}
{\bysame}, 
\textit{Novikov conjecture},
Scholarpedia, 5(10):7912 (2010).

\bibitem{ePcW:85}
{E.K. Pedersen and C. Weibel},
{\it A nonconnective delooping of algebraic K-theory},
in {\it Algebraic and geometric topology}
(A. Ranicki, N. Levitt, and F. Quinn, eds.),
Lecture Notes in Mathematics {\bf 1126},
Springer-Verlag (1985), 166--181.

\bibitem{ePcW:89}
{\bysame}, {\it K-theory homology of spaces}, 
in {\it Algebraic
topology} (G.~Carlsson, R.L.~Cohen, H.R.~Miller, and
D.C.~Ravenel, eds.), Lecture Notes in Mathematics {\bf 1370},
Springer-Verlag (1989), 346--361.

\bibitem{fQ:70}
F. Quinn, 
\textit{A geometric formulation of surgery}, 
in \textit{Topology of Manifolds} (Proc. Inst., Univ. of Georgia,
Athens, Ga., 1969), Markham, Chicago (1970), 500--511.

\bibitem{dRrTgY:11}
{D. Ramras, R. Tessera, and G. Yu},
\textit{Finite decomposition complexity and the integral Novikov conjecture for higher algebraic K-theory},
J. Reine Angew. Math. \textbf{694} (2014), 129--178.

\bibitem{aR:80}
{A. Ranicki},
\textit{Exact sequences in the algebraic theory of surgery},
Princeton U. Press (1980).

\bibitem{aR:98}
{\bysame}, 
{\it Additive L-theory}, K-theory \textbf{3} (1989), 163--195.

\bibitem{aR:92}
{\bysame}, 
{\it Algebraic $L$-theory and topological manifolds}, Cambridge U. Press (1992).

\bibitem{aR:95}
{\bysame}, 
\textit{On the Novikov conjecture}, in
{\it Novikov conjectures, index theory and rigidity} (S.C. Ferry, A. Ranicki, and J. Rosenberg, eds.), 
Cambridge U. Press (1995), 272--337.

\bibitem{aR:01}
{\bysame},
\textit{An introduction to algebraic surgery},
in \textit{Surveys on Surgery Theory, Vol. 2}, Ann. Math. Studies 149, Princeton U. Press (2001), 81--163.

\bibitem{aRmY:95}
A. Ranicki and M. Yamasaki, 
\textit{Controlled K-theory}, 
Topology Appl. \textbf{61} (1995), 1--59.

\bibitem{hRmV:18}
H. Reich and M. Varisco,
\textit{Algebraic K-theory, assembly maps, controlled algebra, and trace methods},
in \textit{Space -- Time -- Matter: Analytic and Geometric Structures}
(J. Br\"{u}ning, M. Staudacher, eds.), 
De Gruyter, Berlin (2018), 1--50.

\bibitem{jR:94}
J. Rosenberg,
\textit{Algebraic K-Theory and its Applications}, 
Graduate Texts in Mathematics 147, Springer-Verlag (1994).

\bibitem{jR:16}
{\bysame},
\textit{Novikov's Conjecture},
in \textit{Open Problems in Mathematics} (J.F. Nash, Jr., M.Th. Rassias, eds.), Springer (2016), 377--402.

\bibitem{cRbS:72}
C.P. Rourke and B.J. Sanderson,
\textit{Introduction to Piecewise-Linear Topology},
Springer-Verlag (1972).

\bibitem{yR:98}
Yu.B. Rudyak,
\textit{On Thom Spectra, Orientability and Cobordism}, Springer-Verlag (1998).

\bibitem{cS:02}
C. Stark,
\textit{Topological rigidity theorems},
in \textit{Handbook of Geometric Topology} (R. Daverman and R. Sher, eds.),
Elsevier (2002), 1045--1084.

\bibitem{rT:83}
{R.W. Thomason},
{\it The homotopy limit problem},
in {\it Proceedings of the Northwestern homotopy theory conference}
(H.R. Miller and S.B. Priddy, eds.),
Cont. Math. {\bf 19} (1983), 407--420.

\bibitem{kV:89}
K. Varadarajan, 
\textit{The finiteness obstruction of C.T.C. Wall},
John Wiley \& Sons (1989).

\bibitem{fW:78}
{F.~Waldhausen}, \textit{Algebraic $K$-theory of generalized free products}, Ann. Math. \textbf{108} (1978), 135--256.

\bibitem{cW:13}
C.A. Weibel,
\textit{The K-book: An Introduction to Algebraic K-theory},
Graduate Studies in Mathematics 145, Amer. Math. Soc., Providence (2013).

\bibitem{sW:90}
{S. Weinberger},
{\it Aspects of the Novikov conjecture},
{\rm Cont. Math.} {\bf 105} (1990), 281--297.

\bibitem{sW:94}
{\bysame}, 
\textit{The Topological Classification of Stratified Spaces},
Chicago Lectures in Mathematics, Chicago U. Press (1994).

\bibitem{mW:02}
M. Weiss,
\textit{Excision and restriction in controlled K-theory},
Forum Math. \textbf{14} (2002), 85--119.

\bibitem{mWbW:01}
M. Weiss and B. Williams,
\textit{Automorphisms of manifolds}, 
in \textit{Surveys on Surgery Theory, Vol. 2}, Ann. Math. Studies 149, Princeton U. Press (2001), 165--220.

\bibitem{jW:77}
J. West, 
\textit{Mapping Hilbert cube manifolds to ANR's: a solution of a conjecture of Borsuk},
Ann. Math. \textbf{106} (1977), 1--18.

\bibitem{gY:98}
G. Yu,
\textit{The Novikov conjecture for groups with finite asymptotic dimension}, Ann. Math. \textbf{147} (1998), 325--355. 

\bibitem{gY:17}
{\bysame}, 
\textit{The algebraic K-theory Novikov conjecture for group algebras over the ring of Schatten class operators}, Adv. Math. \textbf{307} (2017), 727--753.

\end{thebibliography}
\end{document}